\documentclass[reqno,10pt]{amsart}
\usepackage{amsmath,amsthm,amsfonts,color,graphicx}
\usepackage[latin1]{inputenc}
\usepackage[makeroom]{cancel}
\usepackage{tikz}
\usepackage{pgfplots}
\usepackage{subfigure,multirow,float}
\usetikzlibrary{decorations.pathreplacing}
\usepackage{filecontents}
\usepackage{appendix}
\pgfplotsset{compat=1.13}
\usepackage{silence}
\definecolor{lightGray}{RGB}{235,235,235}
\definecolor{orange}{RGB}{255,128,0}
\definecolor{ucib}{RGB}{0,36,105}
\definecolor{mygreen}{RGB}{0,128,0}
\definecolor{lightBlue}{RGB}{102,153,204}

\newtheorem{thm}{Theorem}[section]

\newtheorem{rem}[thm]{Remark}

\DeclareMathAlphabet{\mathpzc}{OT1}{pzc}{m}{it}

\begin{document}
\bibliographystyle{plain}

\title[The SSEM with Chebyshev Polynomials]{The Smooth Selection Embedding Method with Chebyshev Polynomials}

\author{Daniel Agress}
\author{Patrick Guidotti}
\author{Dong Yan}
\address{University of California, Irvine\\
Department of Mathematics\\
340 Rowland Hall\\
Irvine, CA 92697-3875\\ USA }
\email{dagress@uci.edu  and gpatrick@math.uci.edu and dyan6@uci.edu}

\begin{abstract}
We propose an implementation of the Smooth Selection Embedding Method (SSEM), first described in \cite{AG18}, in the setting of Chebyshev polynomials. The SSEM is a hybrid fictitious domain  / collocation method which solves boundary value problems in complex domains by recasting them as constrained optimization problems in a simple encompassing set. Previously, the SSEM was introduced and implemented using a periodic  box (read a torus) using Fourier series; here, it is implemented on a (non-periodic) rectangle using Chebyshev polynomial expansions. This implementation has faster convergence on smaller grids. Numerical experiments will demonstrate that the method provides a simple, robust, efficient, and high order fictitious domain method which can solve problems in complex geometries, with non-constant coefficients, and for general boundary conditions. 
\end{abstract}

\keywords{Fictitious domain methods, embedding methods, numerical solution
    of boundary value problems, boundary value problems as
    optimization problems, high order discretizations of boundary
    value problems, pseudo-spectral methods.}
\subjclass[1991]{}

\maketitle

\section{Introduction}\label{sec:intro}

This paper focuses on the implementation of the SSEM, previously described in \cite{AG18}, in the setting of Chebyshev polynomials. We begin with a brief overview which serves as a motivation and a description of the general SSEM.

\subsection{The SSEM Method}\label{sec:SSEM} We illustrate the SSEM by studying a second order boundary value problem.
\begin{equation}\label{bvp}
  \begin{cases} \mathcal{A} u=f &\text{in }\Omega,\\ \mathcal{B} u=g 
  &\text{on }\Gamma=\partial \Omega. \end{cases}
\end{equation}
Here, $\mathcal{A}$ is a second order differential operator such as, e.g.,  the Laplace operator $-\Delta$, while $\mathcal{B}$ is boundary operator such as, e.g., the trace Dirichlet boundary condition. Efficient and accurate spectral methods exist to solve \ref{bvp} when $\Omega$ is a simple domain, such as a periodic box or a rectangle. However, for a more complicated $\Omega$, these methods are not directly applicable. Fictitious domain methods seek to apply numerical methods to general domains $\Omega$ by embedding into a simpler (fictitious) larger domain, $\mathbb{B}$, for which simple numerical methods  (such as spectal methods, for instance) are available and simple to implement.  One of the fundamental obstacles to such an approach is the fact that the BVP is only defined on $\Omega$, which is a proper subset of $\mathbb{B}$.  Consequently, the original problem only provides an underdetermined set of equations for unknowns defined on the larger set $\mathbb{B}$; in fact, any extension of the solution $u$ of the original BVP is a member of the affine family of solutions to the underdetermined problem. In particular, solutions with low regularity are members of this family; thus, discretizations of the problem will be unable to accurately approximate them. Previous methods have dealt with this issue by smoothly extending the BVP to the entire fictitious domain in such a way that the problem is no longer underdetermined. However, this introduces the difficulty of properly extending all data of the problem while guaranteeing that the original equations are still satisfied, see \cite{Glo94}.

The SSEM, by contrast, treats the original underdetermined problem as a constraint and seeks to find a smooth representative of the affine family of solutions as the minimizer of a constrained optimization problem defined on the whole fictitious domain $\mathbb{B}$. Specifically, given a norm $\| \cdot \|_{\mathcal{S}}$ on $\mathbb{B}$, we seek to solve the constrained optimization problem
\begin{equation}\label{oppb}
  \operatorname{argmin}_{ \left\{\mathcal{C}u= b
  \right\} } \frac{1}{2}\| u\|_{\mathcal{S}} ^2,
\end{equation}
Here, to simplify notation, we have rewritten the entire BVP as one equation 
$$\mathcal{C}u = b, \text{ where } \mathcal{C} = \left( \begin{array}{c}
\mathcal{A} \\ 
\mathcal{B}
\end{array}
\right) \text{ and } b = \left( \begin{array}{c}
f \\ 
g
\end{array}
\right).$$
We emphasize that the operators $\mathcal{A}$ and $\mathcal{B}$ are left in their original form, and only constrain the function on their domains $\Omega$ and $\partial \Omega$, respectively. The norm $\| \cdot \|_{\mathcal{S}}$ on $\mathbb{B}$ is chosen to enforce the desired degree of regularity. In this way, a smooth solution satisfying the constraint $\mathcal{C}u = f$ is selected. Because the selected solution is now smooth, a spectral discretization can approximate it to a high degree of accuracy. As the numerical experiments will demonstrate, the higher the regularity enforced by the norm $\| \cdot \|_{\mathcal{S}}$ on $\mathbb{B}$, the more accurate the approximation will be (compatibly with the expected regularity of the solution, of course.).

\begin{figure}[H]
 \includegraphics[scale=.5]{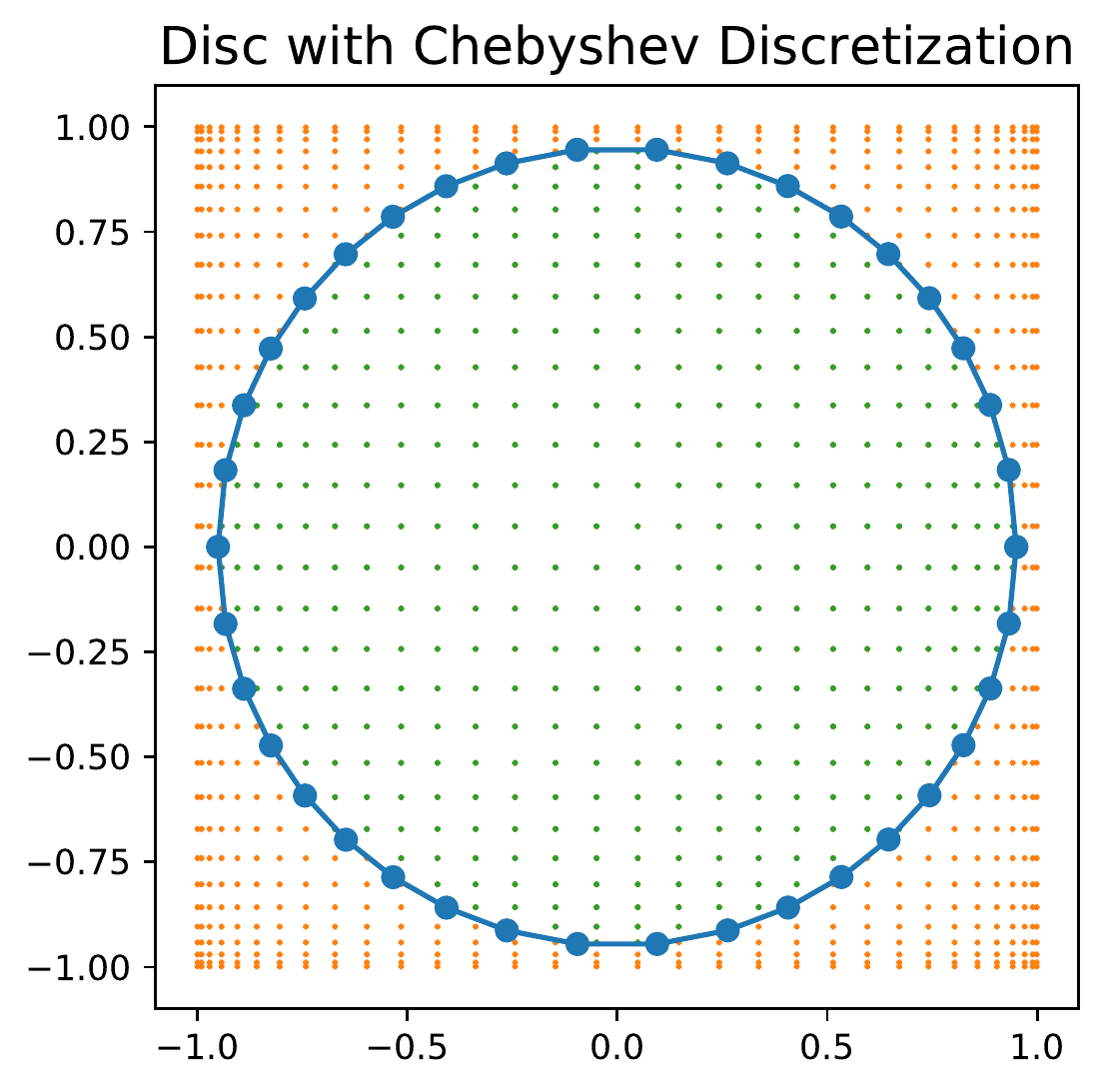}
 \includegraphics[scale=.5]{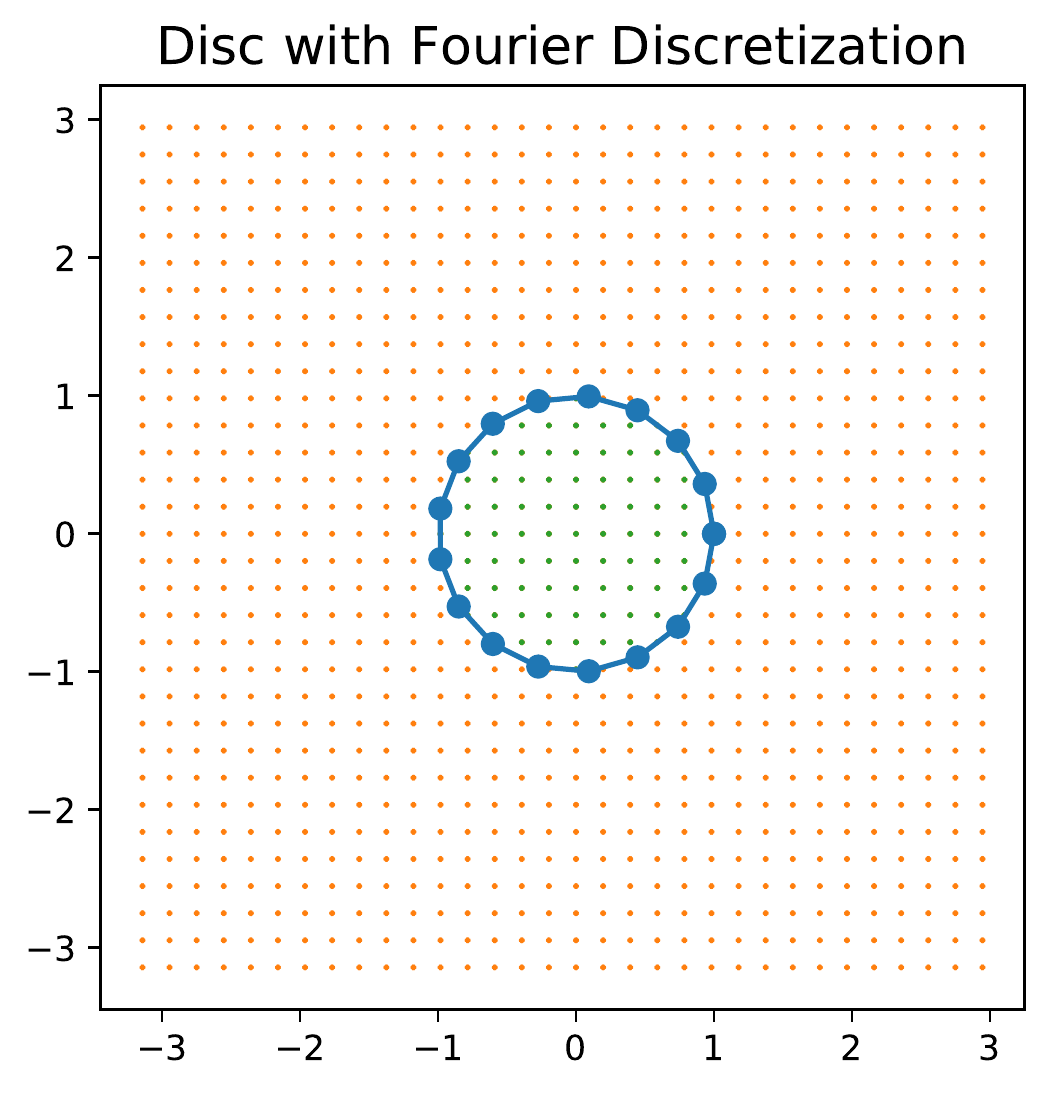}
 \caption{Contrasting the discretization of a disc on the Fourier and Chebyshev grids.}
 \label{fig:domain}
 \end{figure}

We now briefly describe the discretization and the method for solving the optimization problem. To begin, we discretize the encompassing domain $\mathbb{B}$ by a regular grid $\mathbb{B}^m$, selected for the use of spectral methods. Let
$$
\Omega^m = \Omega \cap \mathbb{B}^m
$$
be the discretization of the interior and
$$
\Gamma^m=\{ y_1,\dots,y_{N^\Gamma _m}\}
$$
be a discretization of the boundary obtained by placing roughly equally spaced  points $y_j$ along the boundary. See Figure \ref{fig:domain} for a depiction of two such grids and the corresponding boundary discretizations: one consisiting of Chebyshev roots and one consisting of trigonometric functions' roots, yielding the standard Fourier grid of equally spaced points.
The operators $\mathcal{C}$ and $\mathcal{S}$ are then discretized as
$$
C^m = \begin{bmatrix} A^m\\B^m \end{bmatrix}
$$
and  $S^m$, where $A^m$ and $B^m$ represent an interior  discrete
differential operator at the points of $\Omega^m$ discretizing $\mathcal{A}$ and a discrete realization of the boundary conditions on $\Gamma^m$ encoded by $\mathcal{B}$, respectively. The matrix $S^m$ is a discretization approximating the norm $\mathcal{S}$. Spectral discretizations are chosen for these operators to preserve the accuracy of the method. The original optimization problem can now be described using the discretized operators. Dropping indexes for clarity, we obtain the  problem
\begin{equation}\label{eq:lstsq}
  \operatorname{argmin}_{ \left\{Cu= b
  \right\} } \frac{1}{2}\| u\|_{S} ^2,
\end{equation}
This minimization problem reduces to the following {\em regularized normal equation}
$$
u = S^{-1}C^\top \left( CS^{-1}C^\top \right)^{-1}b,
$$
and the linear system thus obtained can now be efficiently solved in two ways. On the one hand, the matrix $CS^{-1}C^\top$ can be inverted using a preconditioned conjugate gradient method (PCG). This method does not require explicit computation of the matrices and can be used on very dense grids. On the other hand, we observe that, if $(CS^{-1/2})^+$ is the pseudoinverse of $CS^{-1/2}$, the above equation reduces to 
\begin{equation}\label{eq:pinv}
u = S^{-1/2}(CS^{-1/2})^+b.
\end{equation}
The pseudoinverse is efficiently and stably calculated using a QR decomposition of the explicit matrix $CS^{-1/2}$. The advantage of using the pseudoinverse is that the condition number of the matrix $CS^{-1/2}$ is the square root of that of the full matrix $CS^{-1}C^\top$. Thus, higher order regularizing norms can be used, and greater accuracy can be obtained. 

As the numerical experiments performed later will show, the SSEM method can be used to efficiently obtain a spectrally accurate method which works for general geometries, non-constant coefficients, and general boundary conditions. Its implementation, as we shall see, is quite simple as it merely requires discretizing the BVP matrix $\mathcal{C}$ and the regularizing matrix $\mathcal{S}^{-1/2}$ on a regular, rectangular grid. 
\subsection{Motivation for using Chebyshev polynomials.} In \cite{AG18}, the SSEM was implemented for Fourier series on a periodic torus. Such an implementation was chosen because it allowed for the use of the FFT, which made the discretizations straightforward and the computations efficient. While the method was shown to be quite effective at solving BVPs, it suffered from two significant drawbacks. First, because the generated smooth extension of the BVP is periodic, only a small fraction of the periodic domain could be included in $\Omega$. The rest was essentially needed as a buffer in order to allow the extension to smoothly morph into a periodic function. This wasted significant computational resources, because the extension was solved on a far larger grid than was necessary for the solution of the BVP. Put differently, the ratio
$$
\frac{|\Omega^m|}{|\mathbb{B}^m|}
$$
was much smaller than the geometry of $\Omega$ actually required. A second problem was that a $2\pi$-periodic extension of $u$ had far larger optimization norm than $u$ itself. The derivatives outside $\Omega$ were required to be large to force the extension to be periodic. This negatively affected the accuracy of the discretized solution. 

Both of these problems can be remedied by using Chebyshev expansions instead of Fourier series. Instead of working on the periodic torus $[-\pi,\pi)^d$, one can work on the nonperiodic box $\mathbb{B} = [-1,1]^d$. Using Chebyshev polynomials discretized on the Chebyshev roots' grid, functions can still be approximated with spectral accuracy on $\mathbb{B}$. At the same time, because the encompassing domain does not require the periodicity of the extension, far more of the domain $\mathbb{B}$ can be included in $\Omega$. Additionally, for smooth, extendable $u$, the derivative norms of the extension will not be significantly larger than those of $u$ itself. Thus, as demonstrated in the numerical experiments, equivalent accuracy can be achieved on far smaller grids in comparison with the Fourier discretization, leading to faster computation times. As calculations on the Chebyshev grid can also be carried out using the FFT, the method retains the efficiency of the SSEM based on Fourier series. 

\subsection{The Chebyshev smoother}\label{sec:smoother}
The most significant difference between implementing the SSEM on a Chebyshev grid vs on a Fourier grid lies in the choice of regularizing norm used in Equation \ref{oppb}. In the Fourier case, the norm 
$$\| \cdot \|_{S_p} = \| (1-\Delta_{\pi})^{p/2}\cdot \|_{L_2}$$
was used as a smoothing penalty. We note that this norm imposes $H^p$ regularity on the solution, leading, as described in \cite{AG18}, to a $p$ rate of convergence of the error. In addition, the operator $(1-\Delta_{\pi})^{-p/2}$ is diagonalized by the Fourier transform. In particular, the operator $\mathcal{S}^{-1/2}$ has the simple form 
$$ (1-\Delta_{\pi})^{p/2} u = \mathcal{F}^{-1} \circ M \Big[ (1+|k|^2)^{-p/2} \Big] \circ \mathcal{F} u.$$
Here, $(k)_{k \in \mathbb{Z}^d}$ refers to the Fourier frequency vector and $M[\bullet]$ refers to multiplication by the function $\bullet$, where, depending on the occurrence, the function $\bullet$ is either defined everywhere or on the discretization set. The context will clarify which one is meant. The discretized operator $S^{-1/2}$ in \eqref{eq:pinv} is very simply and efficiently computed using the FFT.

Motivated by this choice of smoother for the Fourier case and by the eigenvalue equation described below in Section \ref{sec:eig}, we choose a smoothing norm for the Chebyshev grid accordingly. As described in Section \ref{sec:eig}, define the operator 
$$\mathcal{D} := M \Big[ \sqrt{1-x^2}  \Big] \circ \frac{\partial} {\partial {x}}.$$
Recalling Equation \ref{eq:eig}, we find that for the $m$-th Chebyshev polynomial $T_m(x)$, 
$$ (1-\mathcal{D}^2)^{p/2} T_m(x) = (1+m^2)^{p/2} T_m(x). $$
Exploiting this, we define the norm
$$\| \cdot \|_{\mathcal{S}_p}^2 = \| (1 - \sum_{i=1}^d \mathcal{D}_i^2)^{p/2} \cdot \|^2_{L_2}.$$
Clearly, away from the degeneracies at $-1$ and $1$, this norm imposes $H^p$ regularity on the function $u$. In addition, due to the eigenvalue equation, the operator $\mathcal{S}_p$ is diagonalized by the Chebyshev transform. In particular, if we denote the latter by $\mathfrak{C}$ and let $(k)_{k \in \mathbb{N}^d}$ be the (Chebyshev) frequency vector, we have that
$$
 \mathcal{S}_p^{-1/2}u = \mathfrak{C}^{-1} \circ M \Big[ (1+|k|^2)^{-p/2}\Big] \circ \mathfrak{C} u.
$$
As will be described in Sections \ref{sec:discrete_smoother}, this allows for simple and efficient numerical discretization using the discrete Chebyshev transform as described in \ref{sec:transform}. The numerical experiments of Section \ref{sec:experiments} will demonstrate that, as in the Fourier case, using the $\mathcal{S}_p$ norm leads to a $p$ rate of convergence of the error.
\begin{rem}
Notice that, in the numerical experiments, the observed rate of convergence for the $\mathcal{S}_p$ smoother is somewhat faster than the expected $p$ rate of convergence which was observed in the Fourier case. We suspect that this may have to do with the higher density of points near the boundary of the domain $\Omega$, due to the non-regular spacing of the Chebyshev grid. 
\end{rem}
\begin{rem}
An alternative choice for the norm, which gives a spectral rate of convergence, is
$$
 \| \cdot \|_{S_{\exp}} = \|\mathfrak{C}^{-1} \circ M \Big[ \exp^{\frac{|k|}{2}}\Big] \circ
 \mathfrak{C}u\|_{L_2}.
$$
As motivation, notice that the operator $\mathcal{S}_{\exp}^{-1/2} = \mathfrak{C}^{-1} \circ M \Big[ \exp^{-\frac{|k|}{2}} \Big] \circ \mathfrak{C} u$ is a pseudodifferential operator of heat type. As with $\mathcal{S}_p$, because the operator is diagonalized by the Chebyshev transform, $\mathcal{S}_{\exp}$ can be efficiently discretized and computed.
\end{rem}
\begin{rem}
  In this paper a spectral discretization based on Chebyshev polynomials is chosen. The proposed method can, however, be implemented with respect to any other spectral basis. It is enough to embed $\Omega$ into a larger domain $\mathbb{B}$ for which a full spectral resolution is known for some canonical self-adjoint and positive definite differential operator $D$ with compact resolvent. If the operator admits natural discretizations $\mathbb{B}^m$ for the domain $\mathbb{B}$, $\{ \psi_i \}_{i=1}^m$ for its (orthonormal) eigenfunctions, which are also orthonormal for the appropriate discrete quadrature rule, and satisfy
$$
  \mathcal{D}_m \psi_i^m = \lambda_i^2 \psi_i^m,
$$
for the eigenvalues $\lambda_\bullet^2$ of $D$, then a good smoothing norm given by
$$
 \| \cdot \|_{S_p} = \| \mathfrak{C}_m^{-1} \circ  M \Big[ (1+\lambda_{\bullet}^2)^{p/2} \Big] \circ
 \mathfrak{C}_m \cdot \|_{L_2}
 $$
can be used, where $\mathfrak{C}_m$ is the discrete transformation which computes the coefficients of the (discrete and finite) eigenfunction expansion and $\lambda^2_{\bullet}$ is the corresponding vector of eigenvalues. The Fourier approach of \cite{AG18} clearly fits in this category.
\end{rem}
\subsection{Relationship with the radial basis collocation method}
Even though the SSEM is an example of a so-called embedding method, it can nicely be understood by means of the the so-called radial basis collocation method (RBCM) framework. Surveys of the latter are found in \cite{Fas05} and \cite{Fas07}. In its simplest implementation, the RBCM seeks the solution $u$ to the problem at hand as a linear combination of (smooth) radial functions $\phi$,
$$u = \sum_i c_i \phi(x-x_i),$$
where $u$ is constrained to satisfy the differential operator $\mathcal{C}$ at the points $x_i$. However, in its more general framework, the function $\phi$ is recognized to more properly represent a smoothing kernel applied to a $\delta$ distribution, corresponding to the $\mathcal{S}$ operator in our method, see \cite{Sch07}. In the most general form of the symmetric RBCM method, described in \cite[Chapter~9]{Fas05}, the solution is given by 
$$
u = K \star C^\top \bigl( C(K\star C^\top)\bigr)^{-1}b,
$$
where $K=K(x,y)$ is any symmetric and positive definite smoothing kernel. It is also recognized that such kernels can often be represented as a multiplication operator in Fourier space, leading to both practical and analytical insights for the method. Clearly, the operator $\mathcal{S}$ used in our method fits in the general category of smoothing kernels. In fact, the smoothers $\mathcal{S}_p$ can be viewed as an analog of the Mat\'ern kernels, used in  the RBCM. The Mat\'ern functions (see \cite[Chapter 4]{Fas07}) have Fourier transform $(1+|\xi|^2)^{-p/2}$, while the kernel $\mathcal{S}_p$ has Chebyshev transform $(1+|k|^2)^{-p/2}$. However, to the best of our knowledge, the SSEM distinguishes itself from these other implementations of the RBCM in two main aspects, one practical and one philosophical.
\begin{enumerate}
\item[1.] In collocation methods, the evaluations of the smoothing operator are carried out analytically, and the resulting analytic functions are evaluated at the collocation points. The smoothing kernel itself, thought of as an operator on the encompassing domain $\mathbb{R}^d$, is not explicitly discretized. The SSEM, in contrast, embeds the BVP into a finite domain with a straightforward, regular discretization, e.g. $[-1,1]^d$ with the Chebyshev discretization. The chosen smoothing operators then have simple, discrete approximations on the encompassing domain. Skipping the discretization of the smoother (as in collocation methods) definitely has some advantages. In particular, the collocation points can be placed arbitrarily close, and do not need to take the size of a grid discretization into account. Furthermore, there is no need to perform unnecessary computations on grid points outside of $\Omega$. On the other hand, we believe that having an explicit discretization of the smoother has significant benefits. First, it removes the need for detailed analytic computations of the radial basis functions and their derivatives, and replaces them with a natural multiplication operator in the frequency space of the discretized domain. In fact, for many smoothing operators, a simple explicit representation of the basis functions may not exist. Second, it allows all computations to be done with a straightforward application of the FFT; thus, the evaluation of the matrix requires $O(n\log n)$ as opposed to $O(n^2)$. Finally, having an explicit representation of the smoother allows for the computation of the matrix $CS^{-1/2}$ as opposed to the full matrix $CS^{-1}C^\top$. As described in Section \ref{sec:SSEM}, the condition number of this matrix is the square root of that of the full matrix, which significantly slows the onset of numerical inaccuracy due to ill-conditioning.
\item[2.]Philosophically, as an embedding method, the SSEM looks for the optimal way to embed the solution into an encompassing domain. Collocation methods, in contrast, look for the optimal set of basis functions to place at the collocation points. While in certain cases, the two formulations are equivalent, we believe that the optimization perspective offers many benefits. Firstly, it may often provide natural smoothing operators, suited to the encompassing domain, which might not be apparent at first glance. The norms used in this paper certainly fit in this category. Secondly, the optimization viewpoint allows for the choice of more complex norms which may not be accessible in a collocation framework. For example, non-quadratic objective functionals could be considered. For non-regular problems, weighted norms would be a natural choice. Finally, we note that the optimization perspective is cited in \cite[Chapter~6]{Fas05} as a mathematical justification for the usefulness of the RBCM.
\end{enumerate}

\section{Method}\label{sec:method}
We now detail the implementation of the SSEM on the Chebyshev grid. As described in \eqref{eq:pinv}, the SSEM reduces to solving the equation
$$u = S^{-1/2}(CS^{-1/2})^+ b$$
for $C = \begin{bmatrix} A & B\end{bmatrix}^\top$, a discretization of the BVP and a discretization $S^{-1/2}$ of the smoothing operator described in Section \ref{sec:smoother}. The main step is the description of how to generate the proper discretizations of $C$ and $S^{-1/2}$ and of the solution method used to deal with the resulting linear system.
\subsection{Discretization of the domain.}\label{discretization}
We recall that the BVP is posed in $\Omega \subseteq [-1,1]^d$. The domain $[-1,1]^d$ is discretized by a product set of the Chebyshev grid $\mathbb{B}^m$, described in Section \ref{sec:grid} and given by 
$$
\mathbb{B}^m = \left\{ (x^1, \ldots, x^d) \, \bigg| \, x^i \in \mathbb{C}^m\text{ for } 1\leq i \leq d \right\},
$$
where
$$
 \mathbb{C}^m=\Big\{ \cos \left( \pi \frac{2k+1}{2m} \right) \, \big| \, 0 \leq k \leq m-1 \Big\}.
$$
We can then define the discretized interior as $\Omega^m = \Omega \cap \mathbb{B}^m$, containing $N^m_{\Omega} := |\Omega^m|$ points.  The boundary is discretized by choosing equally spaced points along the boundary $\Gamma$, yielding a set $\Gamma^m$ containing $N^m_{\Gamma} := |\Gamma^m|$ points. We emphasize that these boundary points do not need to lie on the regular grid, but rather lie on the actual boundary $\Gamma$. In two dimensions, the discretization can be achieved by equally spacing points along an arclength parametrization of the boundary curve. In three dimensions, it is not possible to get a perfectly even distribution of points along a two-dimensional boundary surface. However, methods exist to obtain good approximations; see \cite{Pal16} for an example of such an algorithm.

A slightly better boundary discretization, particularly well adapted to the density of the Chebyshev grid, can be obtained as follows. If the boundary $\Gamma$ is a hypersurface contained in $\mathbb{B}$ and parametrized by
$$
\Big(\Gamma_1(\textbf{z}), \ldots, \Gamma_{d}(\textbf{z})\Big), \: \textbf{z} \in S^{d-1},
$$
where $S^{d-1}$ is the $d-1$-dimensional unit sphere, we can create an even distribution of points
$$
 \{ \widetilde{y}\,^i\}_{i=1}^{N^\Gamma_m} = \big\{ (\widetilde{y}\,^i_1,\ldots, \widetilde{y}\,^i_d)\, \big |\, i=1,\dots,N^\Gamma_m \big\}
$$
along
$$
\widetilde{\Gamma} := \Big(\arccos(\Gamma_1(\textbf{z})),\ldots, \arccos(\Gamma_{d}(\textbf{z})) \Big),
$$
now a hypersurface of $[0,\pi)^d$. We note that applying the $\arccos$ function componentwise to the points of $\mathbb{B}^m$ leaves one with a regular grid on $[0,\pi)^d$; thus, setting 
$$
 \Gamma^m = \big\{(y^i_1,\ldots,y^i_d)\big\}_{i=1}^{N^\Gamma_m} = \big\{(\cos(\widetilde{y}\,^i_1),\ldots, \cos(\widetilde{y}\,^i_d))\big\}_{i=1}^{N^\Gamma_m}
$$
produces a boundary discretization of $\Gamma$ the density of which is proportional to the density of the Chebyshev points in $\mathbb{B}$. In the two dimensional numerical experiments, this method was used to discretize the boundary.

A choice is left as to the specific density of the points on the boundary. Increasing the number of discretization points increases the accuracy. However, if the boundary points are more closely spaced than the regular grid points, the regular grid will be unable to distinguish/resolve the boundary points and the matrix will become severely ill-conditioned. In our numerical experiments, we have placed $\frac{m}{2}$ points per unit length in $\widetilde{\Gamma}$ for two dimensional problems. In three dimensions, we have used  $4\pi m^2$ points per unit area on $\Gamma$. These densities seem to provide a good balance of accuracy vs. condition number.
\begin{rem}\label{lobatto}
The grid described above is frequently referred to as the Chebyshev roots grid or the Chebyshev points of the first kind. An alternative choice of grid could have been made with the Chebyshev extrema grid, also known as the Chebyshev points of the second kind. Similar rates of convergence are observed with these points. However, in our numerical experiments, the roots grid appears to be more numerically stable. Furthermore, the regularizer $S_p^{-1}$ described in Section \ref{sec:smoother} is only symmetric for the Chebyshev roots grid, which makes it more convenient for the use of iterative solvers. Henceforth, the Chebyshev grid will refer to the Chebyshev roots grid.
\end{rem}

\subsection{Discretization of the differential operators} \label{sec:matrices}
We recall that the boundary value problem we are solving is of the form
$$
  \begin{cases} \mathcal{A} u=f &\text{in }\Omega,\\ \mathcal{B} u=g 
  &\text{on }\Gamma. \end{cases}
$$
As described in the Introduction, in the SSEM method, the entire BVP, both the interior condition $\mathcal{A}$ and the boundary condition $\mathcal{B}$, are treated as a system of constraints to an optimization problem formulated on the encompassing domain $\mathbb{B}$. This is done by requiring that the discrete solution satisfy an equation corresponding to a discretization $A$ of the operator $\mathcal{A}$ at the points $\Omega^m$, the discretized interior differential equation, and to a discretization $B$ of $\mathcal{B}$ at the points of $\Gamma^m$, the discretized boundary condition. We now describe exactly how to construct these discretizations in the setting of the Chebyshev grid. 
\subsubsection{Constructing $A$}
The matrix $A$ will be an $(N^m_{\Omega} \times m^d)$ matrix which uses values on the entire grid $\mathbb{B}^m$ to approximate the operator $\mathcal{A}$ at the points of $\Omega^m$. For a general second order elliptic BVP, the interior operator is of the form
$$\mathcal{A}u = -a_{ij}u_{x_i x_j} + b_i u_{x_i} + cu.$$
To evaluate the derivatives, we will use the discrete differentiation matrices $D_i$ for the Chebyshev grid described in Section \ref{sec:derivatives}. Here $D_i$ corresponds to differentiation along the $x_i$ direction. Similarly, $D^2_{ij}$ will correspond to taking two derivatives: along the $x_i$ direction and along the $x_j$ direction. We note that, in Section \ref{sec:derivatives}, the differentiation matrices are given implicitly, as linear operators using the discrete cosine and sine transforms $\textsf{DCT}$ and $\textsf{DST}$. This was done to speed the calculations up and to limit RAM usage, although explicit matrices can certainly be used as well. A restriction operator $R:\mathbb{R}^{\mathbb{B}^m} \rightarrow \mathbb{R}^{\Omega^m}$ is also needed which acts by restricting a grid function $u$ to its values on $\Omega^m$. Its transpose, $R^\top : \mathbb{R}^{\Omega^m} \rightarrow \mathbb{R}^{\mathbb{B}^m}$ will act as extension by $0$ from $\Omega^m$ to $\mathbb{B}^m$. For a function $u$ defined on $\mathbb{B}^m$, the operator $A$ is then defined as
$$A(u) = -a_{ij} R \left( D^2_{ij} u \right) + b_i R \left( D_i u \right) + c u.$$
For a function $v$ defined on $\Omega^m$, $A^\top$, used in Section \ref{sec:soln}, is similarly given by
$$A^\top(v) = -D^2_{ij} R^\top \left( a_{ij} v \right) + D_i R^\top \left( b_i v \right) + R^\top \left( c v \right).$$

\begin{rem}\label{location}
While in our implementation, the matrix $A$ enforces the (discrete) differential equation at the points of the Chebyshev grid, this is not strictly necessary. While the skeleton of the method is the Chebyshev grid, the original equations, just like is done for the boundary conditions, can be imposed at any point of $\Omega\subset \mathbb{B}$. Doing so simply requires the use of spectral interpolation beside that of spectral differentiation.
\end{rem}
\subsubsection{Constructing $B$}
The matrix $B$ will be an $(N^m_{\Gamma} \times m^d)$ matrix which uses values on the entire grid $\mathbb{B}^m$ to approximate the operator $\mathcal{B}$ at the points of $\Gamma^m$. For the boundary conditions considered in this paper (Dirichlet, Neumann, or Robin), the boundary operator is of the form
$$
\mathcal{B}u = a \gamma_{\Gamma}u + b \gamma_{\Gamma}(\nabla u)\cdot \nu _\Gamma.
$$
for some smooth functions $a$ and $b$ defined on the boundary $\partial \Omega$. Here, $\gamma_{\partial \Omega}$ is the trace operator and $\nu _\Gamma$ is the unit outward pointing normal vector to $\Gamma$. We note that the choice $a \equiv 1, b \equiv 0$ corresponds to Dirichlet boundary conditions, while $a \equiv 0, b \equiv 1$ corresponds to Neumann boundary conditions. To evaluate the trace and the normal derivative on the boundary, we will need to use the spectral interpolation operators, $\delta_y$ and $\delta_y \circ \nabla$, described in Section \ref{sec:interpolation}. Notice that the vector $\delta_y$ is a spectral discretizations of the $\delta$ distribution located at $y$, while $\nabla$ is the discretized gradient. We construct the matrix $B$ by building each row independently. The $i$-th row of $B$ corresponds to the evaluation of the boundary condition at the $i$-th point of $\Gamma^m$, which we denote as $y_i$. We therefore set
$$[B]_{i\bullet} = a(y_i)\delta_{y_i} + b(y_i)(\delta_{y_i} \circ \nabla) \cdot \nu_{y_i}.$$ 
Here, $\nu_{y_i}$ is the normal vector to $\Gamma$ at the point $y_i$.
\subsubsection{Construction of $S^{-1/2}$}\label{sec:discrete_smoother}
Recall from Section \ref{sec:smoother} that we will be using the smoothers
\begin{align*}
\mathcal{S}_p^{-1/2} &= \mathfrak{C}^{-1} \circ M \Big[ (1+|k|^2)^{-p/2} \Big] \circ \mathfrak{C}
\end{align*}
Here, $\mathfrak{C}$ is the Chebyshev transform, defined in Section \ref{sec:transform} and $(k)_{k\in \mathbb{N}^d}$ is the Chebyshev frequency vector on the $d$ dimensional box $\mathbb{B}$. Each of these can be discretized simply and efficiently using the discrete Chebyshev transform $\mathfrak{C}_m$, also defined in Section \ref{sec:transform}. Defining $(k_{\bullet})_{k_{\bullet} \in \{ 0,\ldots, m-1 \} ^d}$ to be the discrete Chebyshev frequency vector on the $d$ dimensional grid $\mathbb{B}^m$, we define the discrete smoothers 
\begin{align*}
S_p^{-1/2} &= \mathfrak{C}_m^{-1} \circ M \Big[ (1+|k_{\bullet}|^2) ^{-p/2} \Big] \circ \mathfrak{C}_m.
\end{align*} 
\subsection{Solving the resulting system}\label{sec:soln}
Now that we have delineated how to implement the matrices $A$, $B$, and $S_p^{-1/2}$, we turn to solving the linear system. As described in Section \ref{sec:intro}, in the SSEM, the solution $u$ is given by Equation \ref{eq:pinv},
$$ u = S_p^{-1/2} (CS_p^{-1/2})^+ b,$$
where $C = \begin{bmatrix}A&B\end{bmatrix}^\top$ and ${}^+$ denotes the pseudoinverse. As $S_p^{-1/2}$ can be easily computed using the description provided in Section \ref{sec:discrete_smoother}, the remaining issue consists in finding the most efficient and numerically stable method to evaluate $(CS_p^{-1/2})^+b.$ Observe that the matrices $C$ and $S_p^{-1/2}$ are described in the previous section and in the Appendix as linear operators rather than explicit matrices. This would suggest using iterative methods. However, because the matrices are ill-conditioned (see Section \ref{sec:experiments}) we have been unable to find iterative methods which are able to fully capture the accuracy of the method. Instead, in our numerical experiments, we have used the QR decomposition of the explicit matrix $S_p^{-1/2}C^\top$ to find the pseudoinverse. Specifically, if $S_p^{-1/2}C^\top = QR$, then
$$(CS_p^{-1/2})^+ = Q(R^\top)^{-1}.$$
To calculate the explicit matrix $S_p^{-1/2}C^\top$, we simply calculate the $i$-th column of the matrix by evaluating $S_p^{-1/2}C^\top e_i$ using the implicit linear operators described in Section \ref{sec:matrices}. Here, $e_i$ is the $i$-th standard basis vector. We then obtain the solution
$$u = S_p^{-1/2} Q (R^\top)^{-1} b.$$

\begin{rem}
As described in the Section \ref{sec:SSEM}, an alternative method for obtaining the solution is by solving 
$$u = S_p^{-1}C^\top (CS_p^{-1}C^\top)^{-1} b.$$
While the matrix $CS_p^{-1}C^\top$ is more poorly conditioned than $CS_p^{-1/2}$, it has the advantage of being symmetric, and therefore can be inverted using the PCG method. Using the FFT, with $O(n\log n)$ computational complexity, the PCG method works more efficiently and on significantly larger grids than the the direct $QR$ factorization. The drawback is that due to the ill conditioning of the $CS_p^{-1}C^\top$ matrix, only $p \leq 6$ can be used. However, for complicated domains, where a large number of points are necessary to resolve the boundary, the PCG is a more effective method. As a full description of a good preconditioner is given in \cite{AG18}, here, we only offer a brief sketch. We consider the full matrix $CS_p^{-1}C^\top$ as a block matrix, acting on the interior $(\Omega^m)$ and boundary $(\Gamma^m)$ components of $b$ separately.
$$CS_p^{-1}C^\top = \left(\begin{array}{cc}
AS_p^{-1}A^\top & A S_p^{-1}B^\top \\ 
B S_p^{-1}A^\top & B S_p^{-1} B^\top
\end{array} 
\right).$$
The preconditioner $P$ is then a block matrix
$$P = \left(\begin{array}{cc}
C_1^{-1} & 0 \\ 
0 & C_2^{-1}
\end{array} 
\right),$$
where $C_1^{-1}$ is an approximate inverse for $AS_p^{-1}A^\top$ and $C_2^{-1}$ is an approximate inverse for $B S_p^{-1} B^\top$. Because the number of boundary points grows slowly relative to the number of interior points, the matrix $B S_p^{-1} B^\top$ remains small even on denser grids, and it can be inverted directly; thus, $C_2^{-1}$ can just be taken to be $(BS_p^{-1}B^\top)^{-1}$. Next, we note that the operator $AS_p^{-1}A^\top$ is of order $2p-4$; thus, using the preconditioner $(1-\mathcal{D}^2_{\Omega})^{p-2}$ will create a well conditioned operator of order $0$. Here, $\mathcal{D}^2_{\Omega}$ is the operator $\sum_{i=1}^d \mathcal{D}_i^2$ restricted to $\Omega$. We then note that the operator $D^2_{\Omega}$ can be thought of as applying the Laplacian on a uniform grid after composition with the change of variables $\arccos$. Thus, a good discretization of $D^2_{\Omega}$ is obtained by a finite difference approximation of the Laplacian on the grid $\Omega^m$, considering the points $\Omega^m$ as if they were equally spaced. We refer to \cite{AG18} for numerical experiments using the PCG method on a Fourier grid; the results on a Chebyshev grid are similar.
\end{rem}
\section{Numerical Experiments}\label{sec:experiments}
\begin{figure}[H]
\includegraphics[scale=.33]{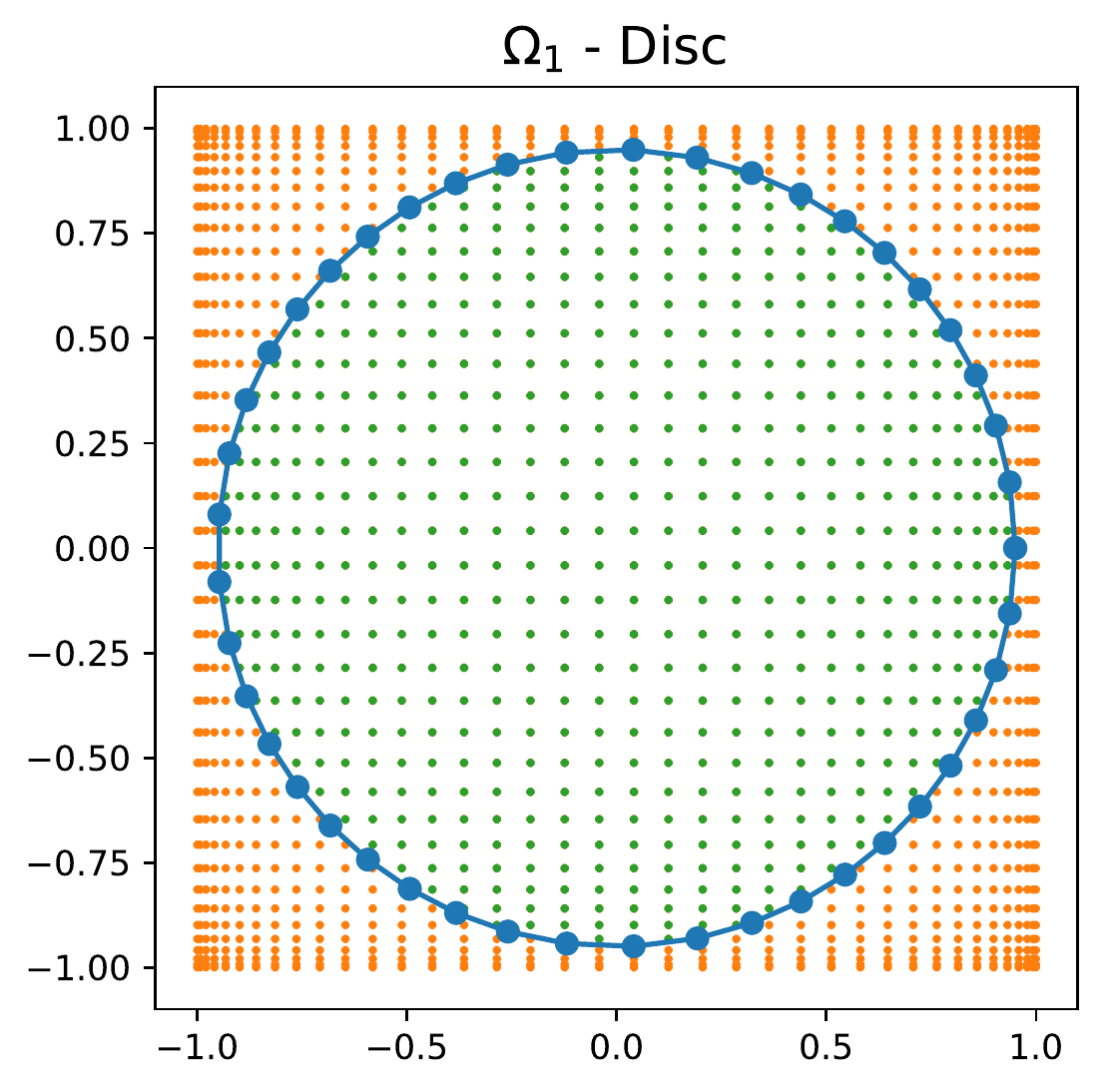}
\includegraphics[scale=.33]{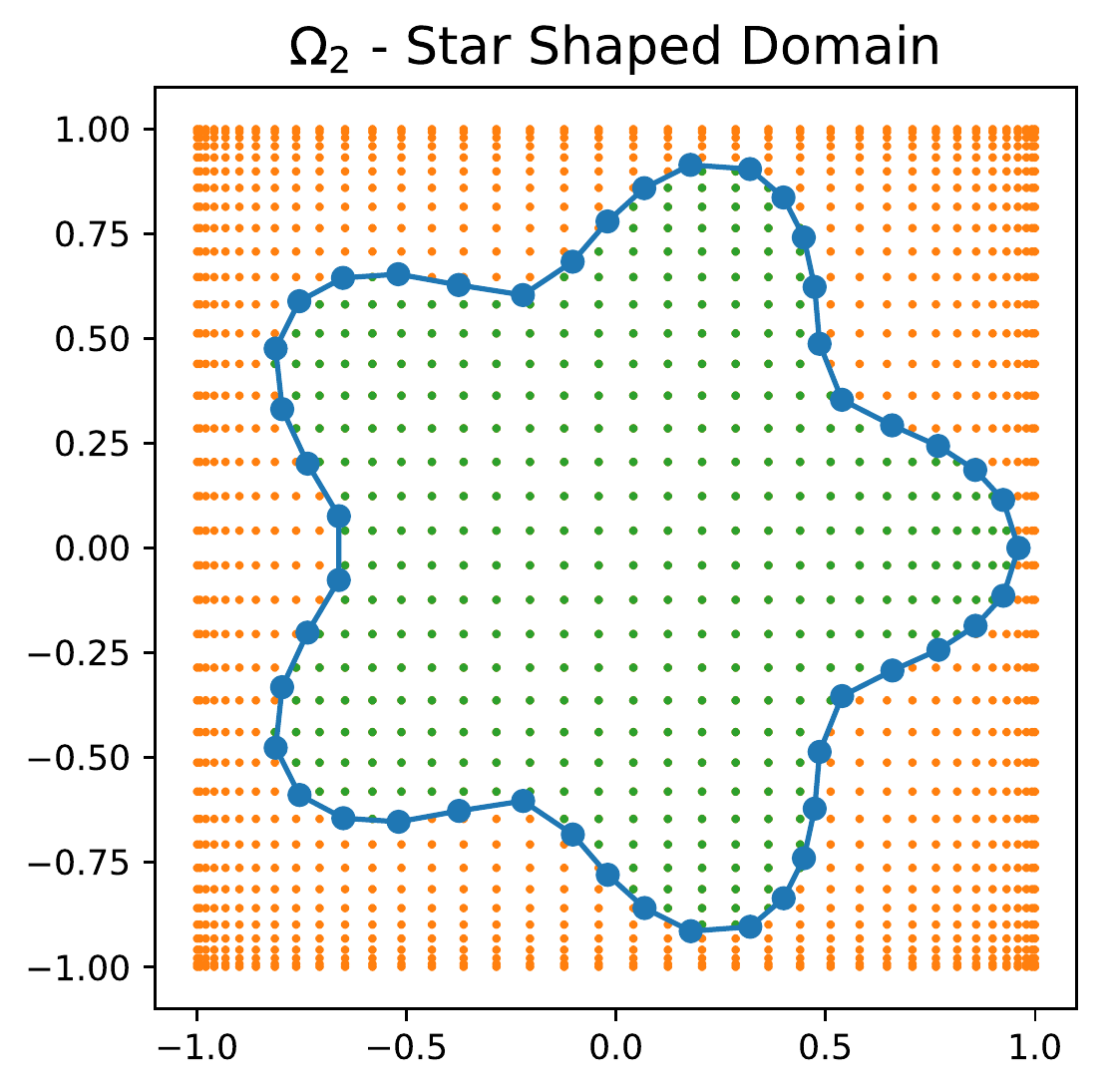}
\includegraphics[scale=.33]{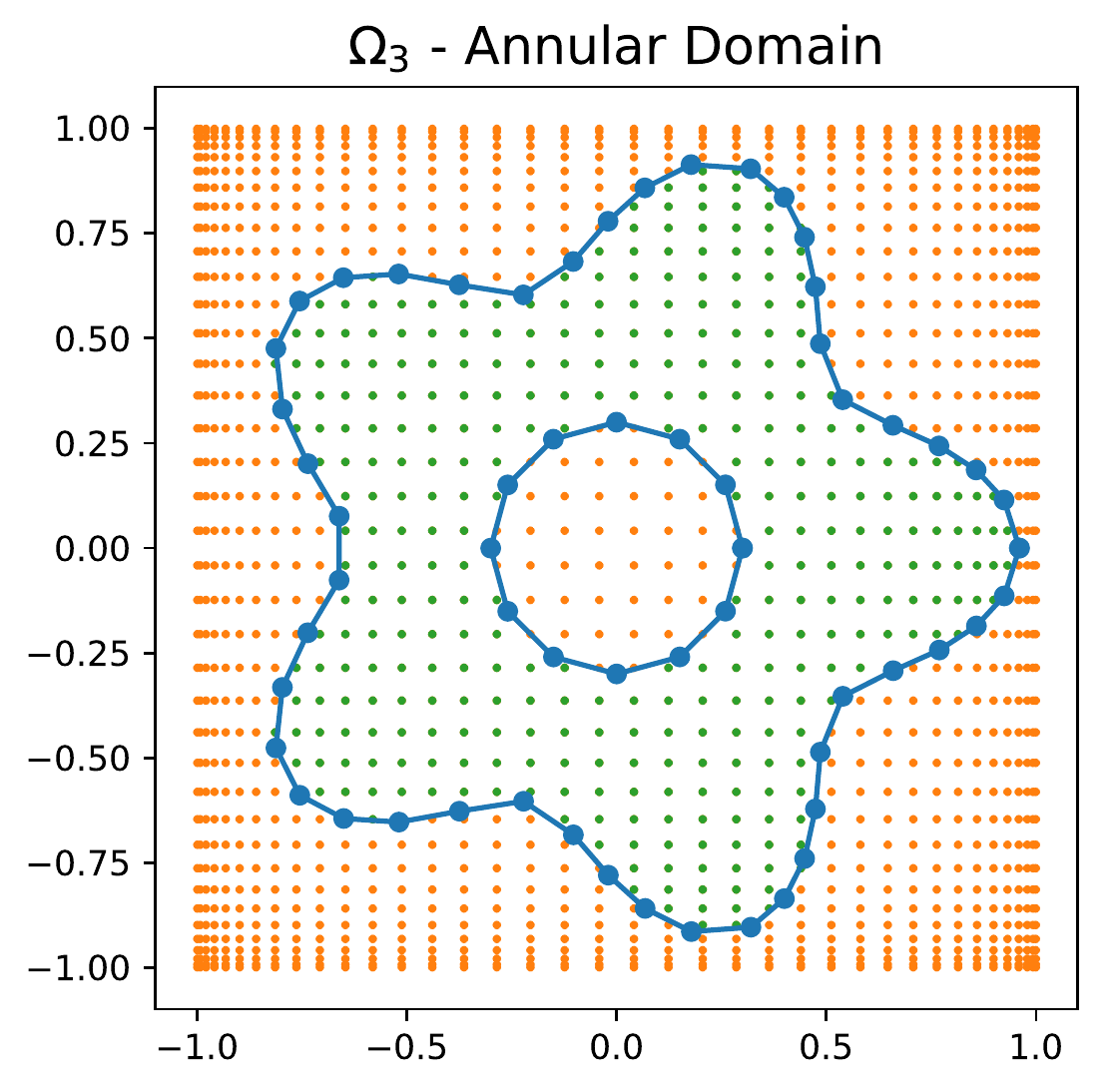}
\caption{Several different domains discretized on the Chebyshev grid considered in the two-dimensional problems.}
\label{fig:2dshapes}
\end{figure}

\begin{figure}[H]
\includegraphics[scale=.33]{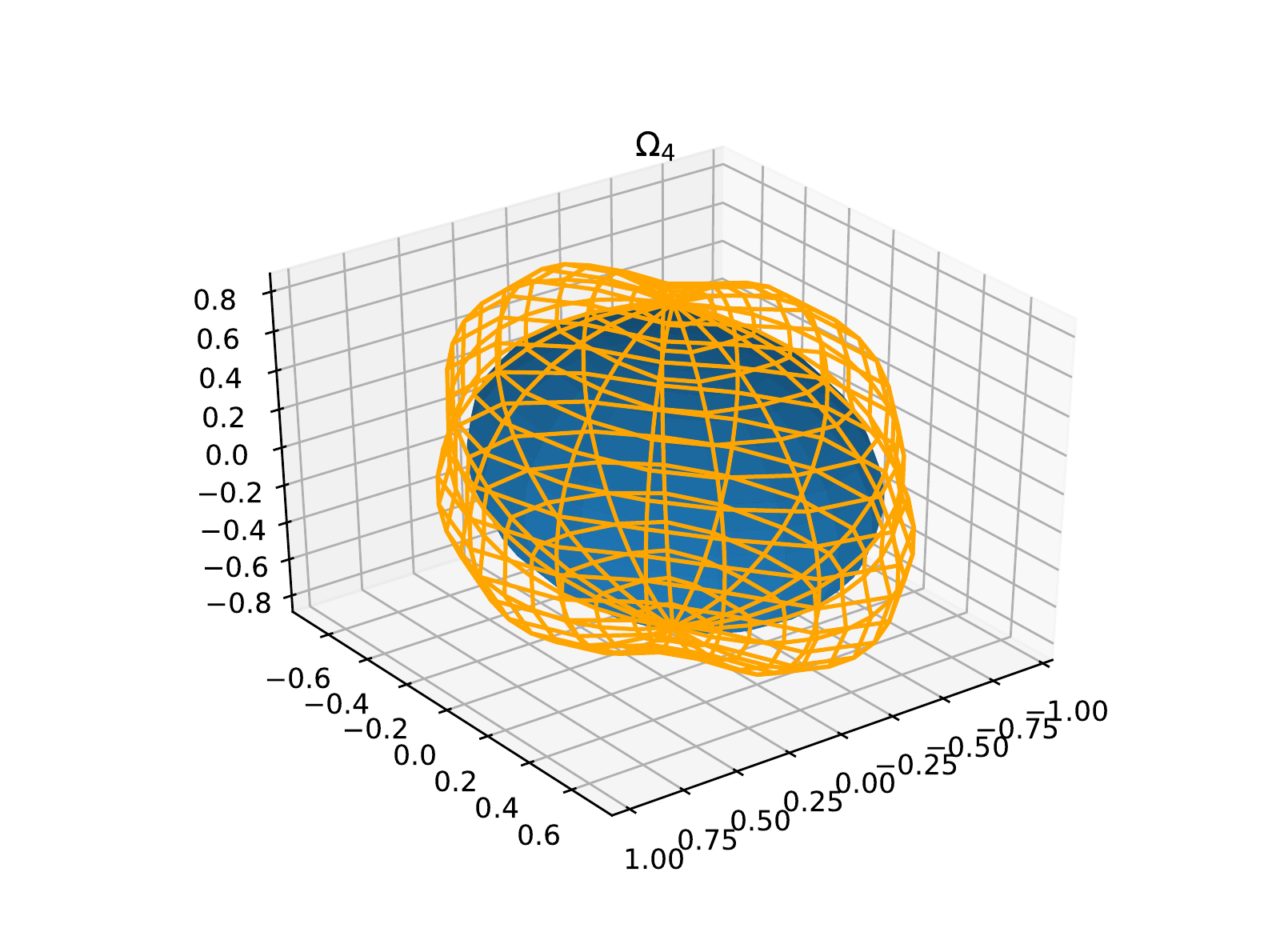}
\includegraphics[scale=.33]{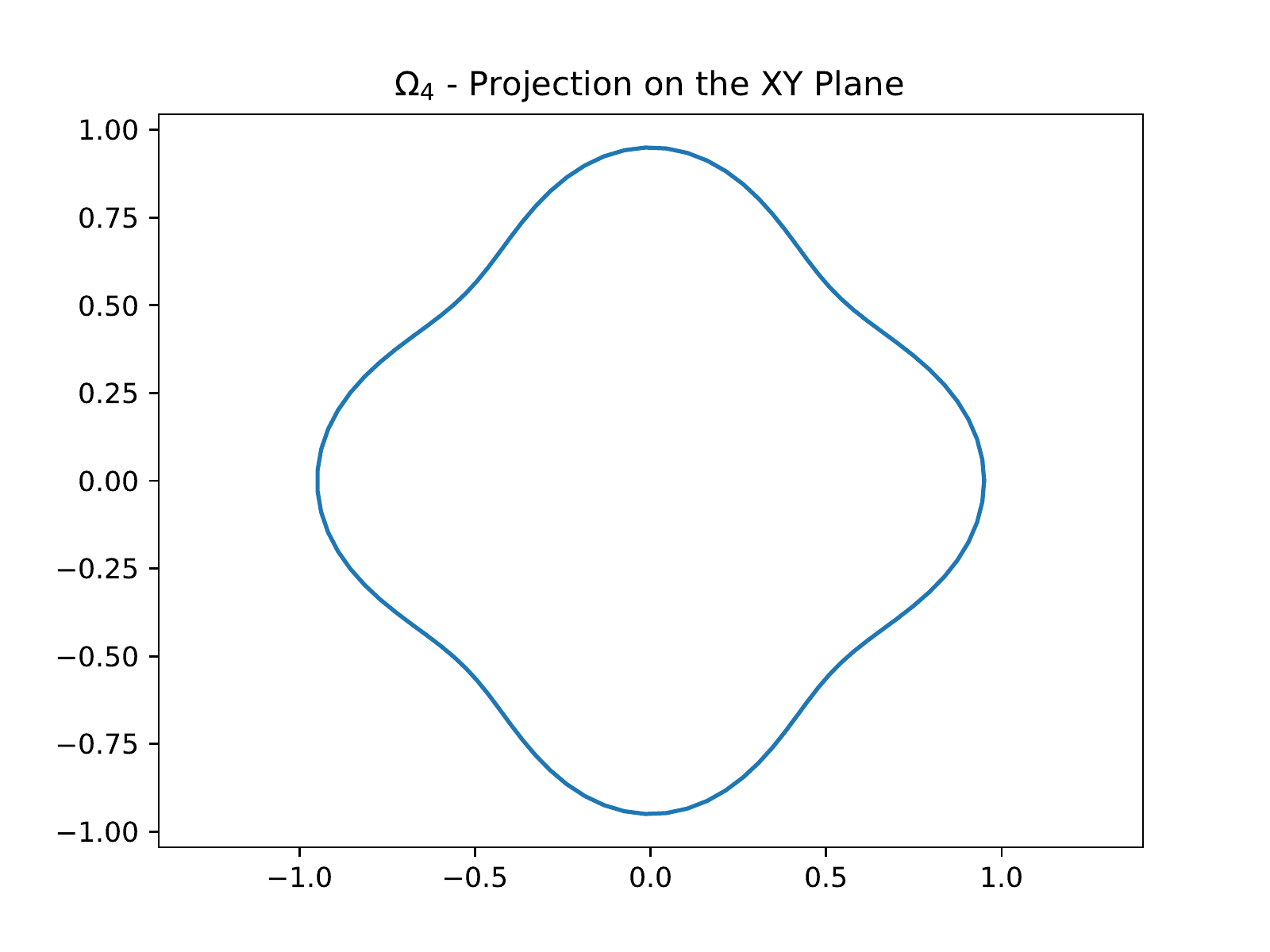}
\caption{Two visualization of $\Omega_4$; a wireframe with a reference sphere on the interior and the projection of $\Omega_4$ onto the XY plane.}
\label{fig:3dshapes}
\end{figure}

We offer a series of numerical experiments to demonstrate the efficacy of the method. We will include Dirichlet, Neumann, and Robin boundary conditions on several different complex two-dimensional domains, and solve a three-dimensional Dirichlet problem. We will also demonstrate how the method can be used to solve a two-dimensional parabolic problem.
For the two-dimensional experiments, we will consider the following domains: a disc, a star-shaped, and an annular domain.
\begin{align*}
 \Omega_1 &= \{(r,\theta) \; \big| \; r < .95 \}.\\
 \Omega_2 &= \{(r,\theta) \; \big| \; r < .8(1+.2\cos(\theta)) \}.\\
 \Omega_3 &= \{ (r, \theta) \; \big| \; .3 < r < .8(1+.2\cos(\theta)) \}.
\end{align*}
For the three-dimensional experiment, we will consider the star-shaped perturbation of the sphere given by
$$\Omega_4 = \{(r,\theta,\phi) \; \big| \; r < .85 + .1\sin(\phi)\cos^2(\theta) \}.$$
The domains are shown in Figures \ref{fig:2dshapes} and \ref{fig:3dshapes}.
\subsection{A two-dimensional Dirichlet problem}\label{sec:dirichlet}

\begin{figure}[H]
\includegraphics[scale=.45]{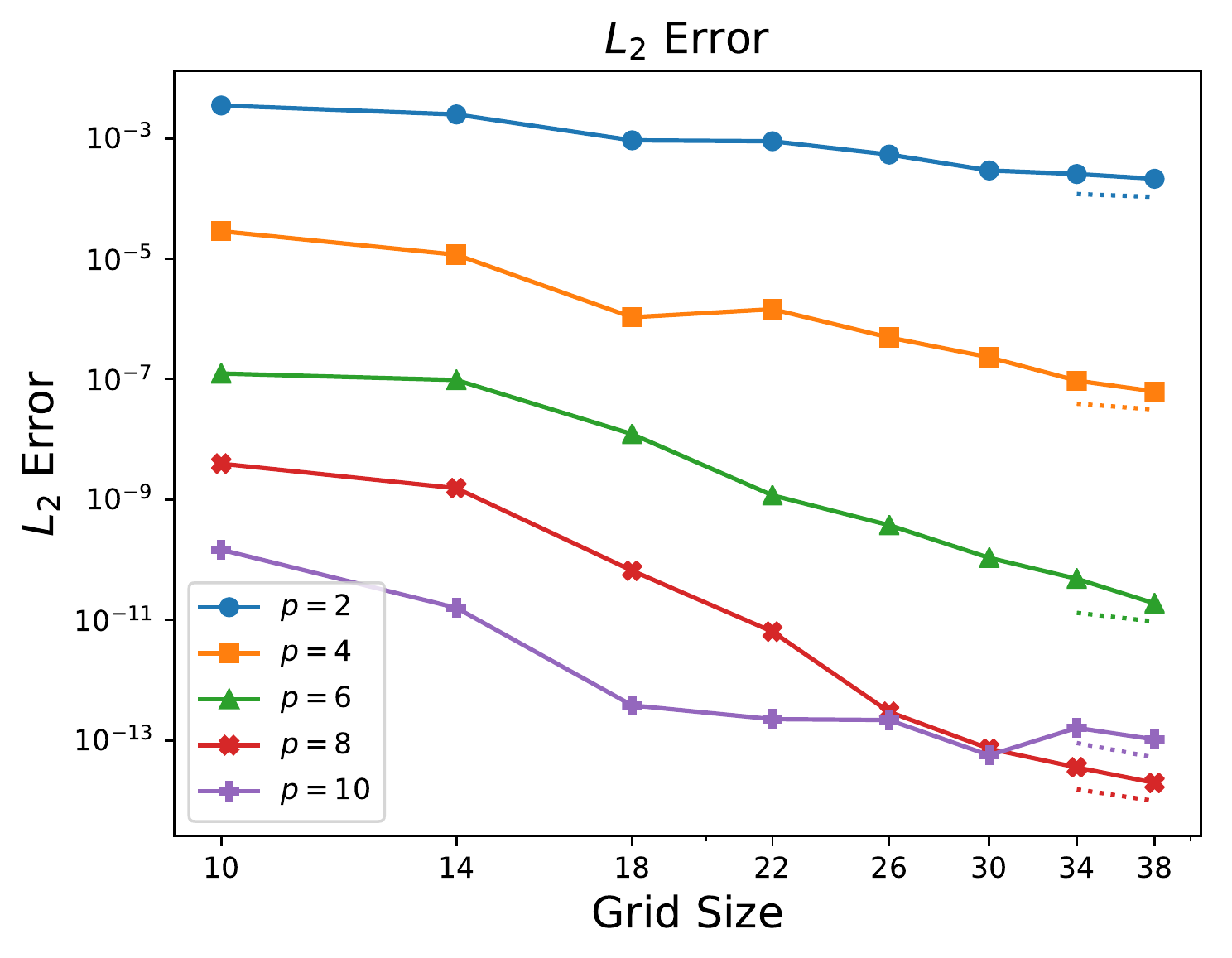}
\includegraphics[scale=.45]{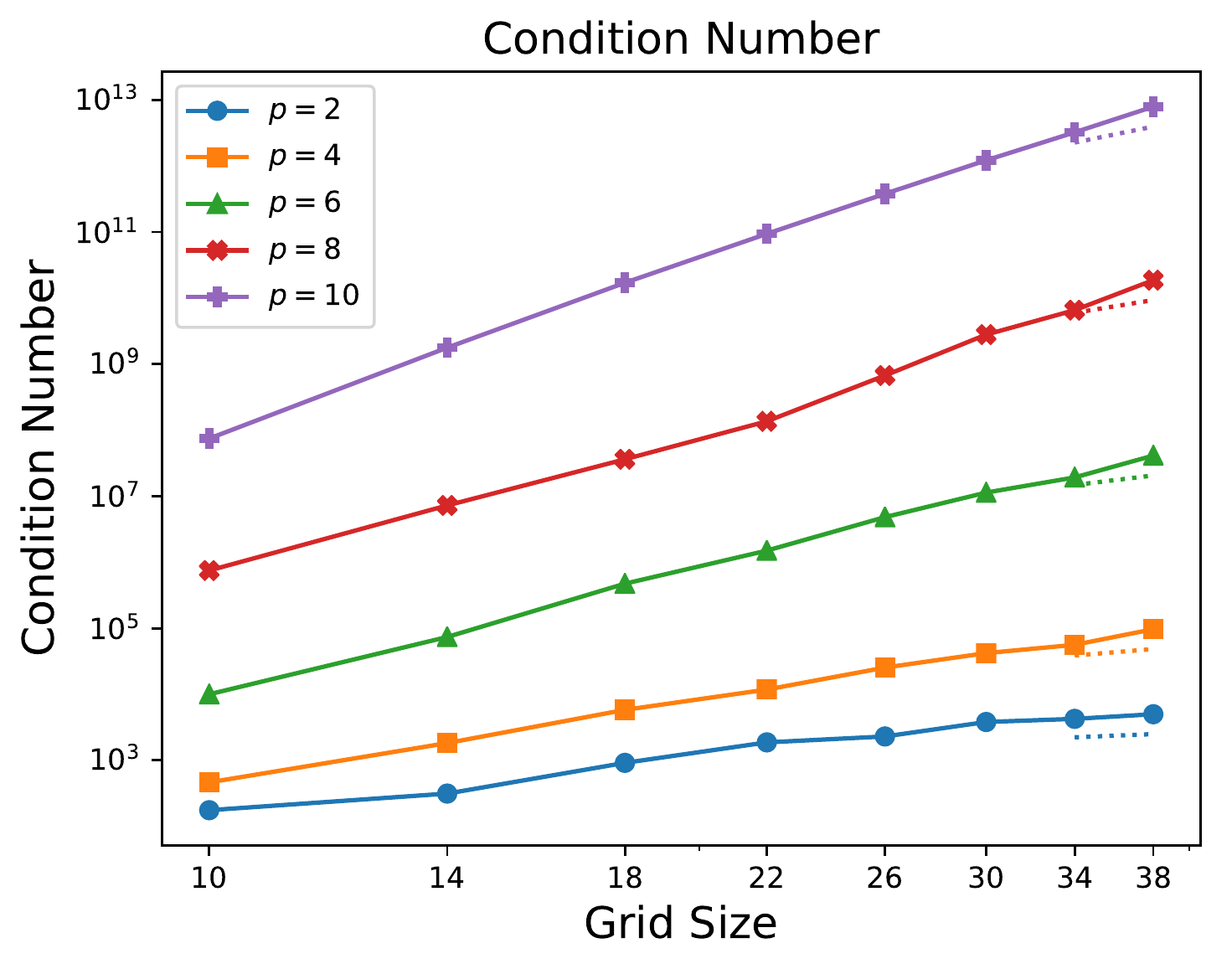}
\caption{Convergence of the $L_2$ error and 
growth of the condition number of the matrix $CS_p^{-1/2}$ for different values of $p$ when solving Equation \protect\eqref{eq:expdiri}. The light dotted lines are reference lines of slope $\frac{1}{m^p}$ and $m^p$, where $m$ is the number of grid points along one dimension.}
\label{fig:ratesdiri}
\end{figure}

On $\Omega_1$, the disc of radius $0.95$ shown in Figure \ref{fig:2dshapes}, we consider the Dirichlet problem
\begin{equation}\label{eq:expdiri}
\begin{cases}
-\Delta u  = 0 \quad  & {\rm in } \,\,\Omega_1,\\
u = x^2 - y^2 \quad  &{\rm on}\,\, \partial \Omega_1.
\end{cases}
\end{equation}
The exact solution is $x^2-y^2$. The sizes of the different discretizations considered, together with the CPU times, are listed in Table \ref{tab:timesdiri}. We note that because we are using explicit matrices, the CPU time is independent of the chosen smoother. In Figure \ref{fig:ratesdiri}, we show the condition number of the resulting system and the convergence rates of the $L_2$ error for the different smoothers $S_p$.
\begin{table}[H]
\centering
\begin{tabular}{|c|c|c|c|}
\hline
$|\mathbb{B}^m|$ & $N^m_{\Omega}$ & $N^m_{\Gamma}$ & CPU Time\\
\hline 
$ 10 ^2 = 100$ & $ 40 $ & $ 12 $ & $ 0.2 $\\
\hline
$ 14 ^2 = 196$ & $ 76 $ & $ 17 $ & $ 0.2 $\\
\hline
$ 18 ^2 = 324$ & $ 136 $ & $ 22 $ & $ 0.2 $\\
\hline
$ 22 ^2 = 484$ & $ 204 $ & $ 26 $ & $ 0.3 $\\
\hline
$ 26 ^2 = 676$ & $ 280 $ & $ 31 $ & $ 0.4 $\\
\hline
$ 30 ^2 = 900$ & $ 372 $ & $ 36 $ & $ 0.5 $\\
\hline
$ 34 ^2 = 1156$ & $ 464 $ & $ 40 $ & $ 0.7 $\\
\hline
$ 38 ^2 = 1444$ & $ 580 $ & $ 45 $ & $ 1.1 $\\
\hline
\end{tabular}
\caption{Grid sizes and computation times for solving Equation \protect\eqref{eq:expdiri}. All computations were performed on an Intel I7-6660U.}
\label{tab:timesdiri}
\end{table}
\subsection{A two-dimensional Neumann problem}\label{sec:neumann}
\begin{figure}[H]
\includegraphics[scale=.45]{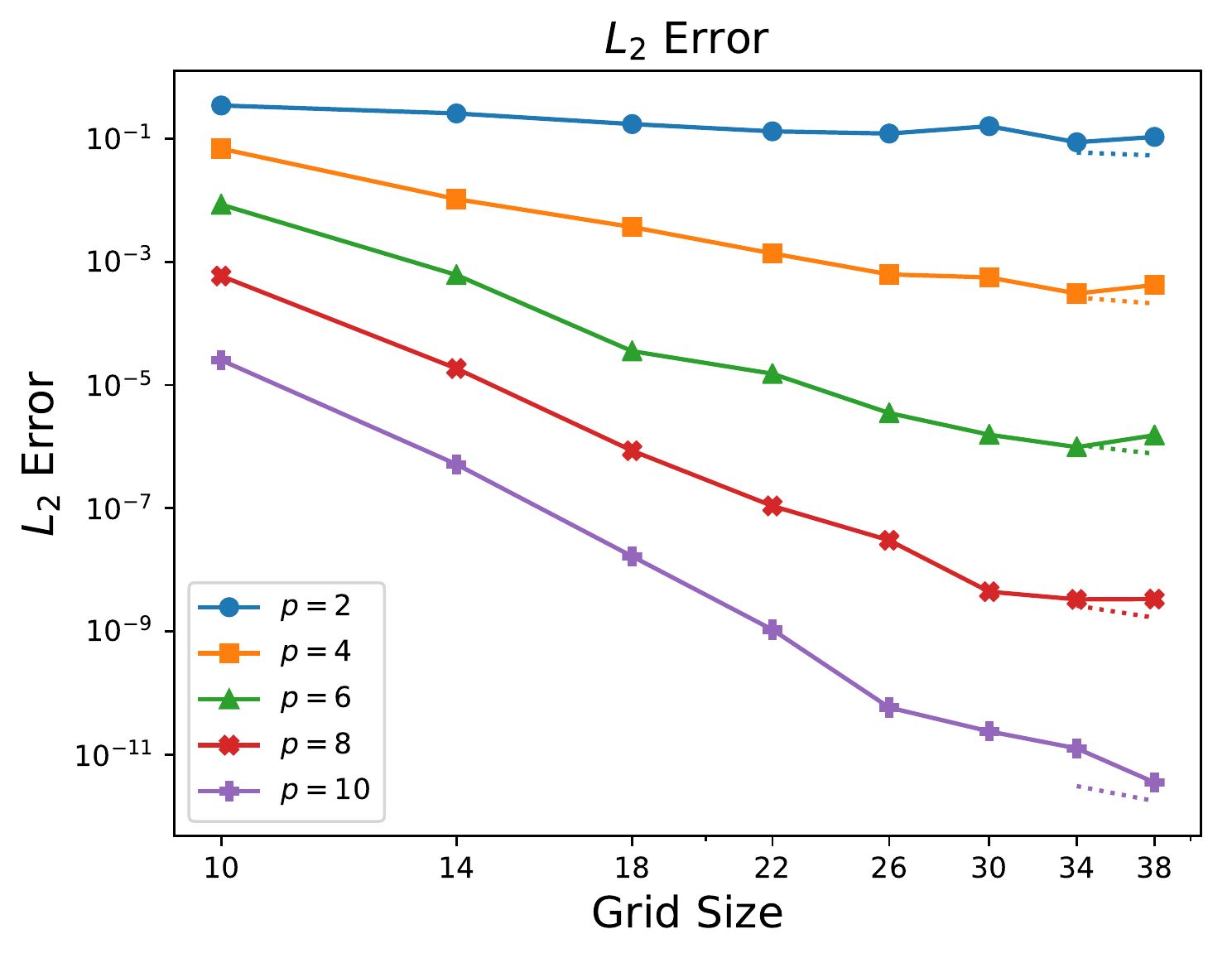}
\includegraphics[scale=.45]{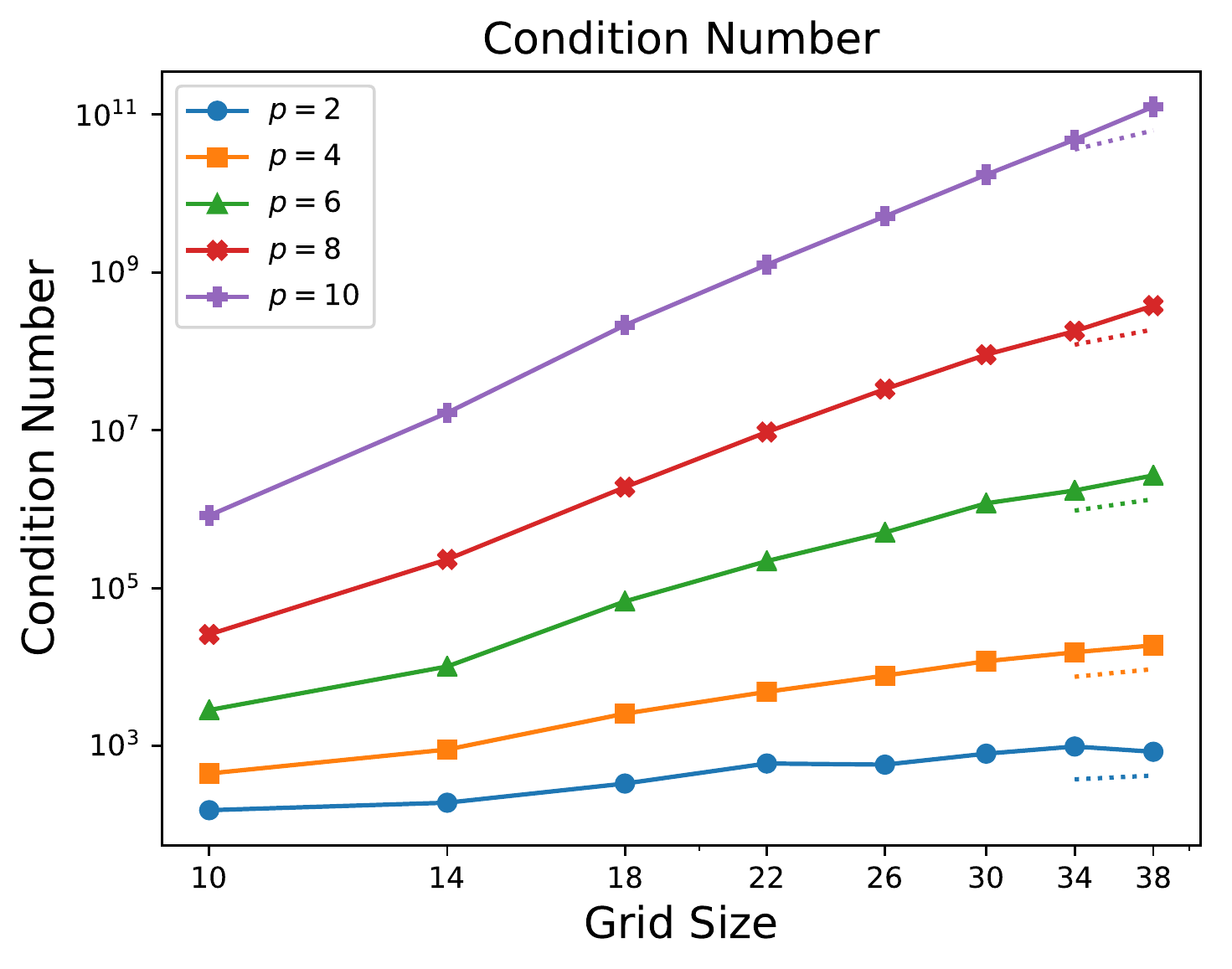}
\caption{Convergence of the $L_2$ error and growth of the condition number of the matrix $CS_p^{-1/2}$ for different values of $p$ when solving Equation \protect\eqref{eq:expneu}. The light dotted lines are reference lines of slope $\frac{1}{m^p}$ and $m^p$, where $m$ is the number of grid points along one dimension.}
\label{fig:ratesneu}
\end{figure}

We now consider the star-shaped domain
$$\Omega_2 = \{(r,\theta)\mid r<.8(1+.2\cos(5\theta))\}.$$
The discretization is shown in Figure \ref{fig:2dshapes}. We consider the Neumann problem 
\begin{equation}\label{eq:expneu}
\begin{cases}
-\left[(2+y)\partial_x^2 + (2-x)\partial_y^2\right] u  = -12(x+y) \quad  &{\rm in } \,\, \Omega_2,\\
\frac{\partial u}{\partial \nu} = \langle 3x^2,3y^2\rangle \cdot \nu  \quad  &{\rm on}\,\, \partial \Omega_2,
\end{cases}
\end{equation}
where $\nu$ is the outward pointing unit normal vector to $\partial \Omega_2$. The exact solution is $u(x,y) = x^3+y^3$. The numerical results are summarized in Table \ref{tab:timesneu} and Figure \ref{fig:ratesneu}.
\begin{table}[H]
\centering
\begin{tabular}{|c|c|c|c|}
\hline
$|\mathbb{B}^m|$ & $N^m_{\Omega}$ & $N^m_{\Gamma}$ & CPU Time\\
\hline 
$ 10 ^2 = 100$ & $ 28 $ & $ 13 $ & $ 0.2 $\\
\hline
$ 14 ^2 = 196$ & $ 50 $ & $ 18 $ & $ 0.2 $\\
\hline
$ 18 ^2 = 324$ & $ 86 $ & $ 23 $ & $ 0.2 $\\
\hline
$ 22 ^2 = 484$ & $ 134 $ & $ 28 $ & $ 0.3 $\\
\hline
$ 26 ^2 = 676$ & $ 182 $ & $ 33 $ & $ 0.4 $\\
\hline
$ 30 ^2 = 900$ & $ 242 $ & $ 37 $ & $ 0.4 $\\
\hline
$ 34 ^2 = 1156$ & $ 308 $ & $ 42 $ & $ 0.6 $\\
\hline
$ 38 ^2 = 1444$ & $ 384 $ & $ 47 $ & $ 0.9 $\\
\hline
\end{tabular}
\caption{Grid sizes and computation times for solving Equation \protect\eqref{eq:expneu}. All computations were performed on an Intel I7-6660U.}
\label{tab:timesneu}
\end{table}
\subsection{A two-dimensional Robin problem}
\begin{figure}[H]
\includegraphics[scale=.45]{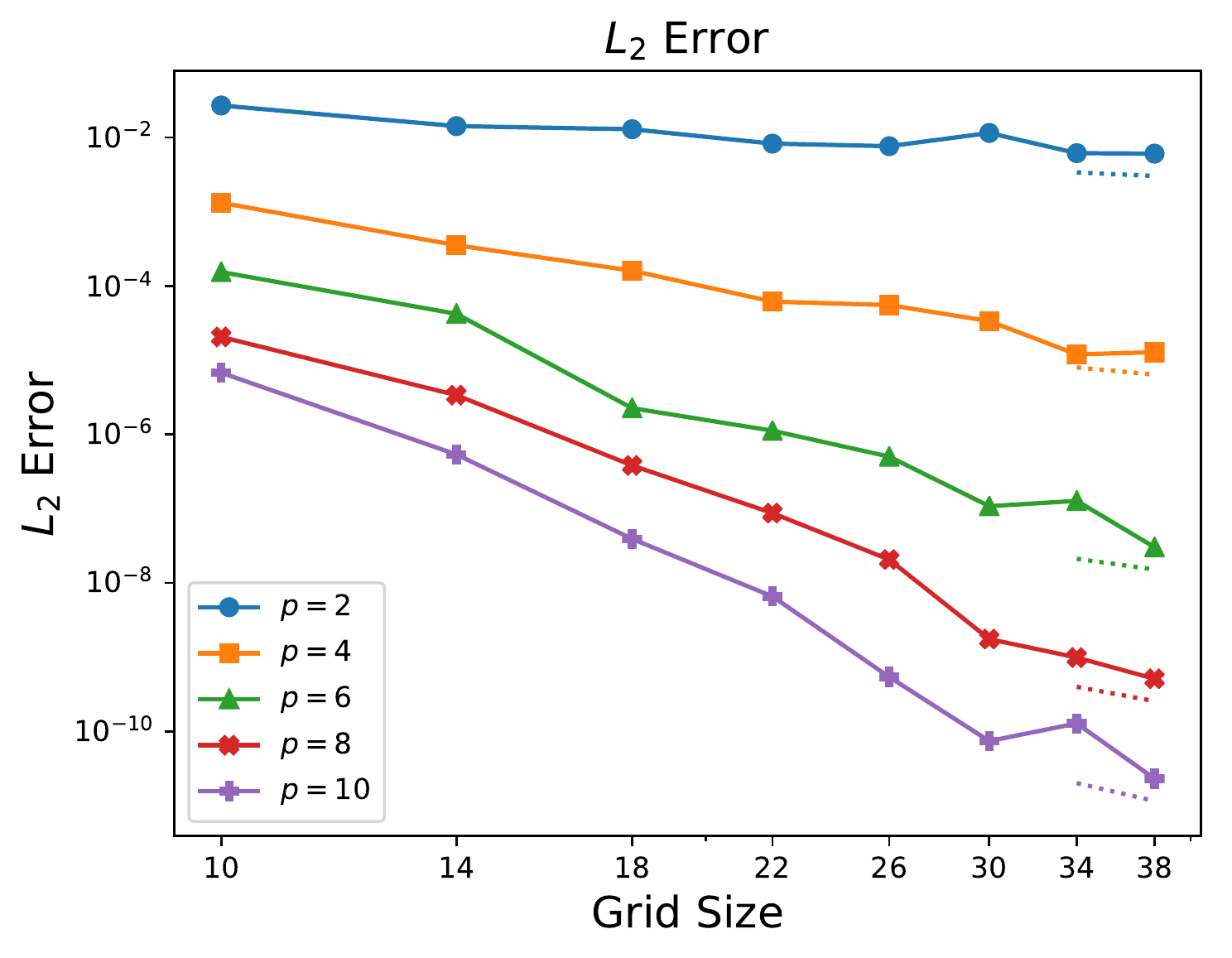}
\includegraphics[scale=.45]{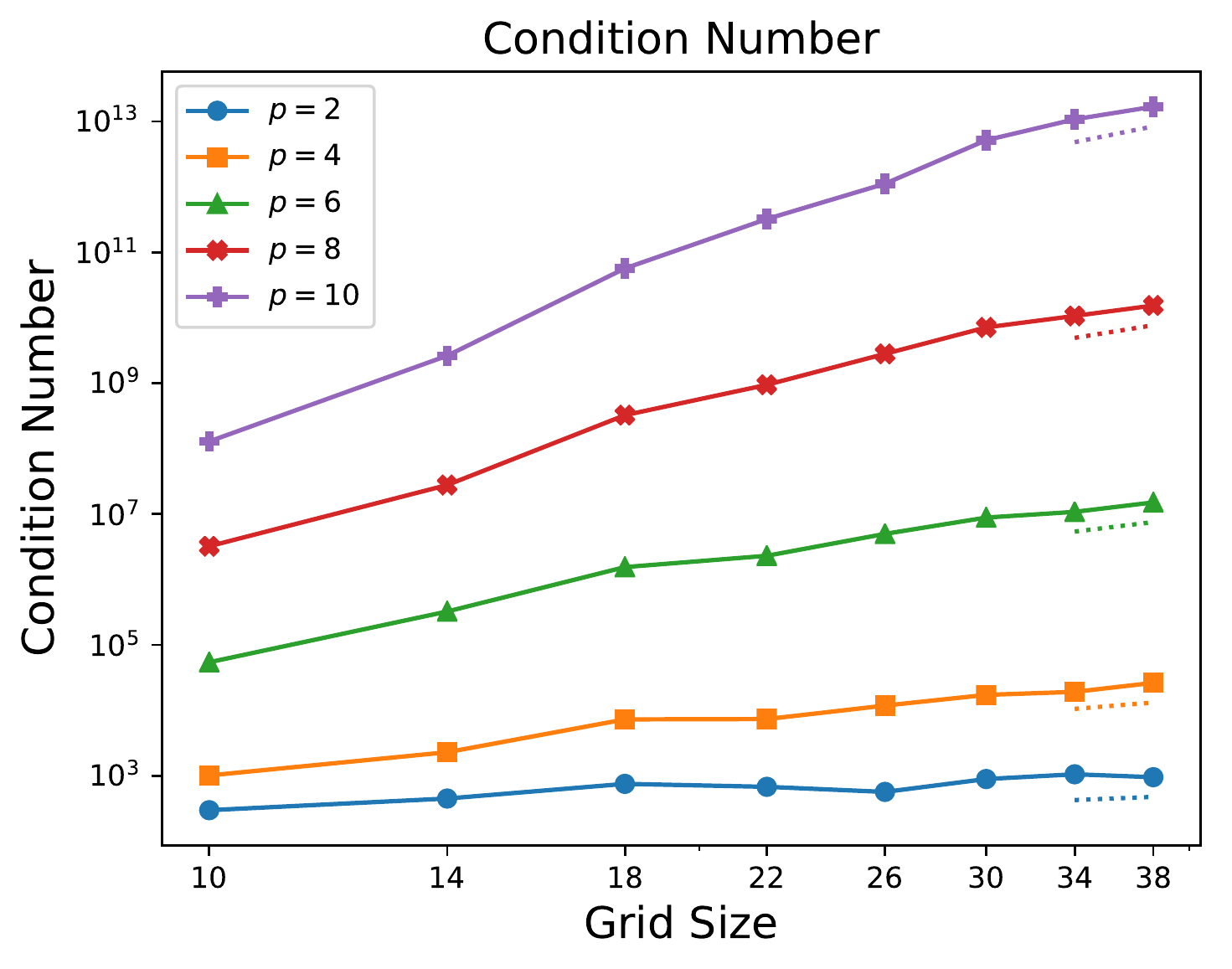}
\caption{Convergence of the $L_2$ error and growth of the condition number of the matrix $CS_p^{-1/2}$ for different values of $p$ when solving Equation \protect\eqref{eq:exprob}. The light dotted lines are reference lines of slope $\frac{1}{m^p}$ and $m^p$, where $m$ is the number of grid points along one dimension.}
\label{fig:ratesrob}
\end{figure}

Next consider the annular domain
$$\Omega_3 = \{(r,\theta) \; \big| \; .3 < r < .8(1+.2(\cos(5\theta))\},$$
shown in Figure \ref{fig:2dshapes}. Consider the Robin boundary value problem
\begin{equation}\label{eq:exprob}
\begin{cases}
-\Delta u  = -\sinh x-\cosh y \quad  &{\rm in } \,\, \Omega_3,\\
u + \frac{\partial u}{\partial \nu} = \sinh x+\cosh y+\langle \cosh x,\sinh y\rangle \cdot \nu \quad  &{\rm on}\,\, \partial \Omega_3,
\end{cases}
\end{equation}
where $\nu$ is the outward pointing unit normal vector to $\partial \Omega_3$. The exact solution is given by $u(x,y) = \sinh x + \cosh y$. The numerical results are summarized in Table \ref{tab:timesrob} and Figure \ref{fig:ratesrob}.
\begin{table}[H]
\centering
\begin{tabular}{|c|c|c|c|}
\hline
$|\mathbb{B}^m|$ & $N^m_{\Omega}$ & $N^m_{\Gamma}$ & CPU Time\\
\hline 
$ 10 ^2 = 100$ & $ 24 $ & $ 17 $ & $ 0.3 $\\
\hline
$ 14 ^2 = 196$ & $ 46 $ & $ 23 $ & $ 0.3 $\\
\hline
$ 18 ^2 = 324$ & $ 74 $ & $ 29 $ & $ 0.4 $\\
\hline
$ 22 ^2 = 484$ & $ 122 $ & $ 35 $ & $ 0.4 $\\
\hline
$ 26 ^2 = 676$ & $ 166 $ & $ 41 $ & $ 0.5 $\\
\hline
$ 30 ^2 = 900$ & $ 218 $ & $ 47 $ & $ 0.7 $\\
\hline
$ 34 ^2 = 1156$ & $ 276 $ & $ 53 $ & $ 0.9 $\\
\hline
$ 38 ^2 = 1444$ & $ 340 $ & $ 59 $ & $ 1.1 $\\
\hline
\end{tabular}
\caption{Grid sizes and computation times for solving Equation \protect\eqref{eq:exprob}. All computations were performed on an Intel I7-6660U.}
\label{tab:timesrob}
\end{table}
\subsection{A three dimensional Dirichlet problem}
\begin{figure}[H]
\includegraphics[scale=.45]{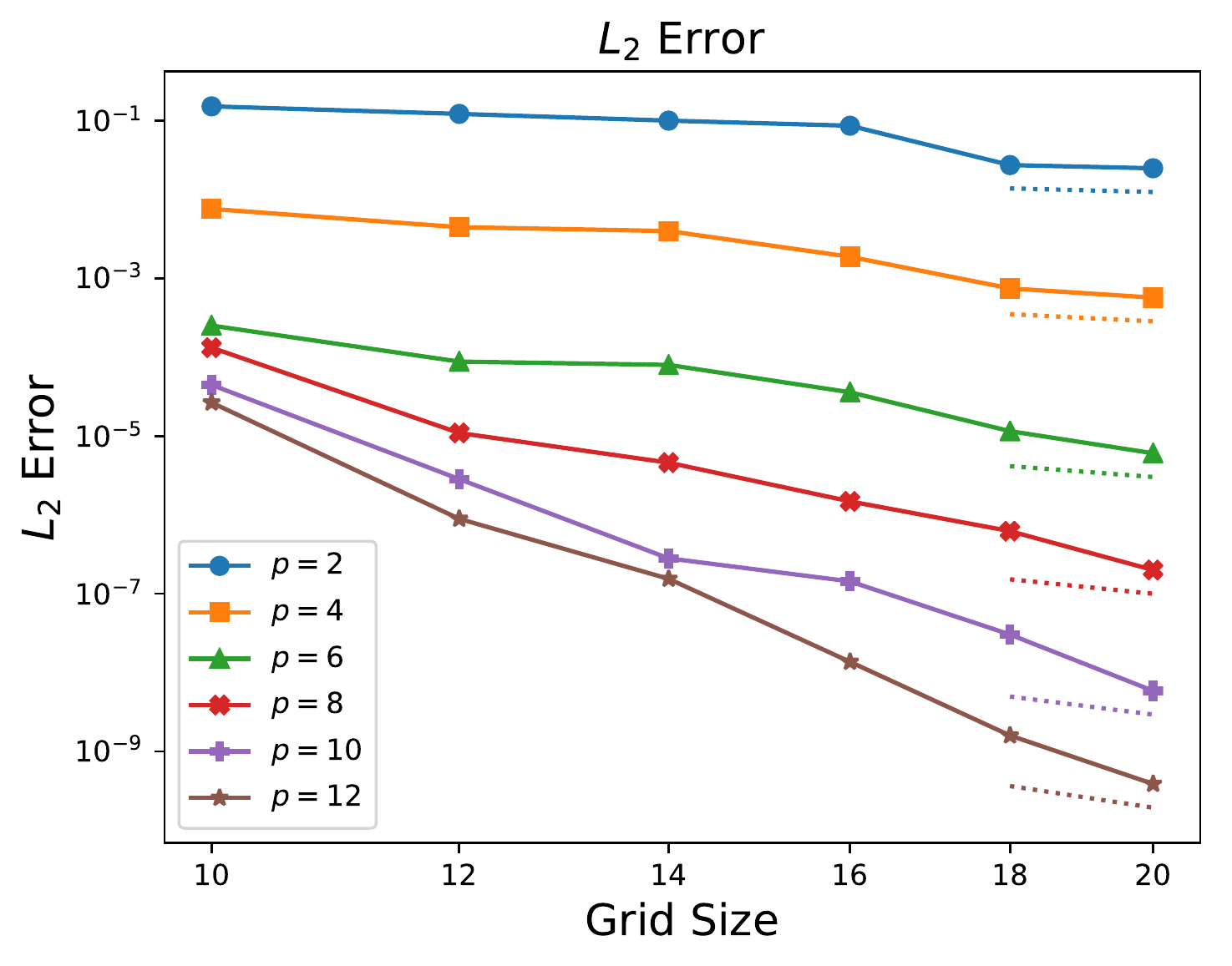}
\includegraphics[scale=.45]{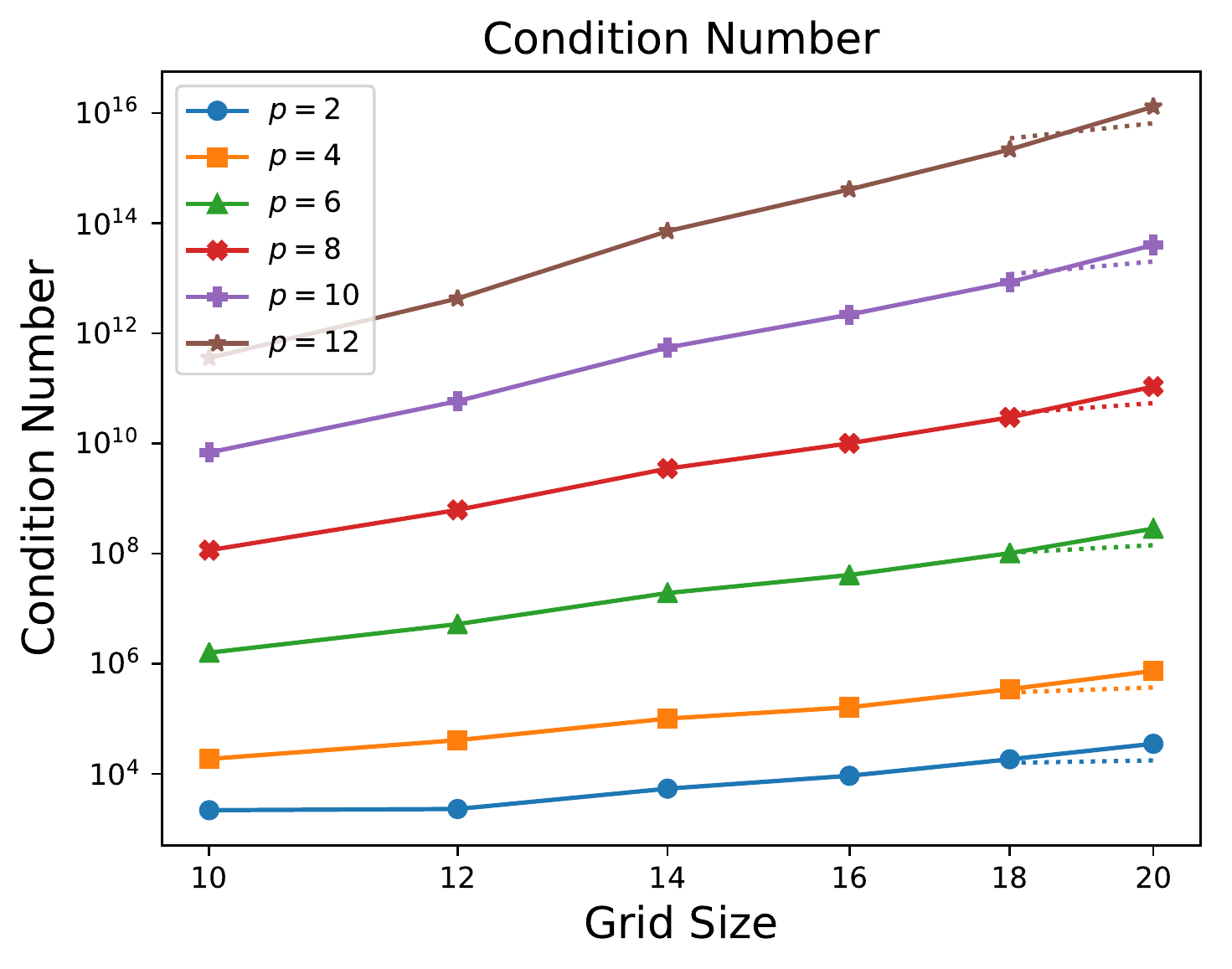}
\caption{Convergence of the $L_2$ error and growth of the condition number of the matrix $CS_p^{-1/2}$ for different values of $p$ when solving Equation \protect\eqref{eq:exp3d}. The light dotted lines are reference lines of slope $\frac{1}{m^p}$ and $m^p$, where $m$ is the number of grid points along one dimension.}
\label{fig:rates3d}
\end{figure}

We consider the star-shaped domain
$$\Omega_4 = \{ (r,\phi,\theta) \; \big| \; r < .85 + .1\sin(\phi)\cos(4\theta) \}.$$
As the domain is a perturbation of the sphere, we have discretized the boundary $\partial \Omega_4$ by projecting a Fibonacci lattice of the unit sphere (see \cite{Gon10}) onto $\partial \Omega_4$. As the Fibonacci lattice is evenly distributed over the sphere, this gives a good discretization of the surface $\partial \Omega_4$. We study the Dirichlet problem
\begin{equation}\label{eq:exp3d}
\begin{cases}
-\Delta u  = -2y\sin(x) + x^2y\sin(z) \quad  & {\rm in } \,\,\Omega_4,\\
u = x^2y\sin(z) \quad  &{\rm on}\,\, \partial \Omega_4.
\end{cases}
\end{equation}
The exact solution is $u(x,y,z) = x^2y\sin(z)$. The numerical results are summarized in Table \ref{tab:times3d} and Figure \ref{fig:rates3d}.
\begin{table}[H]
\centering
\begin{tabular}{|c|c|c|c|}
\hline
$|\mathbb{B}^m|$ & $N^m_{\Omega}$ & $N^m_{\Gamma}$ & CPU Time\\
\hline 
$ 10 ^3 = 1000$ & $ 128 $ & $ 90 $ & $ 0.4 $\\
\hline
$ 12 ^3 = 1728$ & $ 192 $ & $ 129 $ & $ 0.7 $\\
\hline
$ 14 ^3 = 2744$ & $ 296 $ & $ 176 $ & $ 1.7 $\\
\hline
$ 16 ^3 = 4096$ & $ 472 $ & $ 231 $ & $ 4.5 $\\
\hline
$ 18 ^3 = 5832$ & $ 672 $ & $ 292 $ & $ 9.3 $\\
\hline
$ 20 ^3 = 8000$ & $ 896 $ & $ 361 $ & $ 17.4 $\\
\hline
\end{tabular}
\caption{Grid sizes and computation times for solving Equation \protect\eqref{eq:exp3d}. All computations were performed on an Intel i7-7700HQ.}
\label{tab:times3d}
\end{table}

\subsection{A parabolic problem}
Finally we describe a method to solve a time dependent problem using the SSEM. We consider the parabolic cylinder $\Omega \times [0,2]$, where $\Omega$ is the star-shaped domain described in Section \ref{sec:neumann} and pictured in Figure \ref{fig:2dshapes}. Define $j_0(r)$ to be the $0$-th Bessel function of the first kind. Letting $r$ denote the Euclidean distance from $0$, consider the radial function
$$u(r,t) =  e^{-t} j_0(r) - e^{-\frac{t}{4}} j_0\Big(\frac{r}{2}\Big).$$
The function $u$ then satisfies the parabolic BVP
\begin{equation}\label{eq:exppara}
\begin{cases}
u_t - \Delta u = 0 & \text{in } \Omega\times (0,2],\\
u = j_0(r) - j_0(\frac{r}{2}) & \text{in } \Omega \times \{0\}, \\
u = e^{-t}j_0(r) - e^{-\frac{t}{4}}j_0(\frac{r}{2}) & \text{on } \Gamma \times (0,2].
\end{cases}
\end{equation}
To discretize the domain, we use the Chebyshev grid $\mathbb{B}^m$ described in Section \ref{sec:grid}. As for the time interval $[0,2]$, we use a (shifted) Chebyshev extrema grid, 
$$
\mathbb{B}^{n}_E = \{ t_j \}_{j=0}^{n}\text{, where } t_j = -\cos \Big(\frac{\pi j}{n}\Big) + 1.
$$
In this section $n=10$ was chosen in all of the experiments. The full discretization of the parabolic cylinder $[-1,1]^2 \times [0,2]$ is then given by $\mathbb{B}^m \times \mathbb{B}^n_E$.  The choice to use the extrema grid rather than the standard Chebyshev (roots) grid in the time variable was made because imposing the boundary condition at $t=0$ is slightly more straightforward, since the boundary point $t=0$ lies on the grid. For a description of how to construct the time differentiation matrix, $D_t$, we refer to \cite{Tre00}. 

With the discretization $\mathbb{B}^m \times \mathbb{B}^{n}_E$, the interior of the parabolic cylinder is given by $\Omega^m \times \widetilde{\mathbb{B}}^n_E$, where
$$
\widetilde{\mathbb{B}}^n_E = \big\{t_i \in \mathbb{B}^n_E \, \big| \, t_i > 0\big\}.
$$
The discretized ``bottom'' boundary of the cylinder is given by $\Omega^m \times \{0\}$, whereas the discretization of the lateral boundary $\Gamma \times (0,2]$ is simply given by $\Gamma^m \times \widetilde{\mathbb{B}}^n_E$. Letting $R_{K^m}$ denote the evaluation operator on the discrete set $K^m$, we can define the matrices
\begin{align*}
A &= R_{\Omega^m \times \widetilde{\mathbb{B}}^n_E} \circ (D_t - D_{x_1}^2 - D_{x_2}^2),\\
B_1 &= R_{\Omega^m \times \{0\}},\\
B_2 &= R_{\Gamma^m \times \widetilde{\mathbb{B}}^n_E}.
\end{align*}
Notice that evaluation of a function on the above sets simply amounts to their restriction to the sets since $\Omega^m \times \widetilde{\mathbb{B}}^n_E$ and $\Omega^m \times \{0\}$ are sets of regular grid points. However, because $\Gamma^m$ does not contain regular grid points in general, the evaluation matrix $R_{\Gamma^m \times \widetilde{\mathbb{B}}^n_E}$ will require the use of the interpolation operators described in Sections \ref{sec:matrices} and \ref{sec:interpolation}. If $b_1$ and $b_2$ represent the evaluations of the function $e^{-t}j_0(r) - e^{-\frac{t}{4}} j_0\Big( \frac{r}{2} \Big)$ at the points of $\Omega^m \times \{0\}$ and $\Gamma^m \times \widetilde{\mathbb{B}}^n_E$, respectively, the BVP is fully discretized by the matrix equation $Cu=b$ where
$$
C = \begin{bmatrix} A & B_1 & B_2 \end{bmatrix}^\top\text{ and }b
= \begin{bmatrix} 0 & b_1 & b_2 \end{bmatrix}^\top
$$

As for the elliptic SSEM, the problem is converted to a constrained optimization problem
$$  \operatorname{argmin}_{ \left\{Cu= b
  \right\} } \frac{1}{2}\| u\|_{S} ^2,$$
where $\| \cdot \|_S$ is a smoothing norm. In the parabolic case, $\| \cdot \|_S$ will need to be a space-time norm over $\mathbb{B}^m \times \mathbb{B}^n_E$. We recall from Section \ref{sec:smoother} and \ref{sec:eig} the operator
$$
(\mathcal{D}^m)^2 =  \mathfrak{C}_m^{-1} \circ M \Big[|k_{\bullet}|^2 \Big] \circ\mathfrak{C}_m.
$$
The operators $\mathcal{D}^m_i$ and $\mathcal{D}^m_t$ will represent applying the operator in the $x_i$ and $t$ directions, respectively. Motivated by our choice of smoothing norm used in the elliptic case and described in Section \ref{sec:smoother}, the norm given by
$$
 \| u \|_{S_p} =\Big\| \Big(1 - \sum_{i=1}^2 (\mathcal{D}^m_i)^2 - (\mathcal{D}^m_t)^2\Big)^{p/2}(u)\Big\|,
$$
is used in order to enforce space-time regularity of the numerical solution. As with the norms described in the elliptic case, this norm has the benefit of simple implementation using the discrete Chebyshev transform. We remark that while this norm is clearly effective, as demonstrated by our numerical experiments, it is not natural from the point of view of parabolic PDEs and may not be the optimal one to use; we are continuing to investigate the best choice of smoother in the parabolic case.

\begin{figure}[H]
\includegraphics[scale=.45]{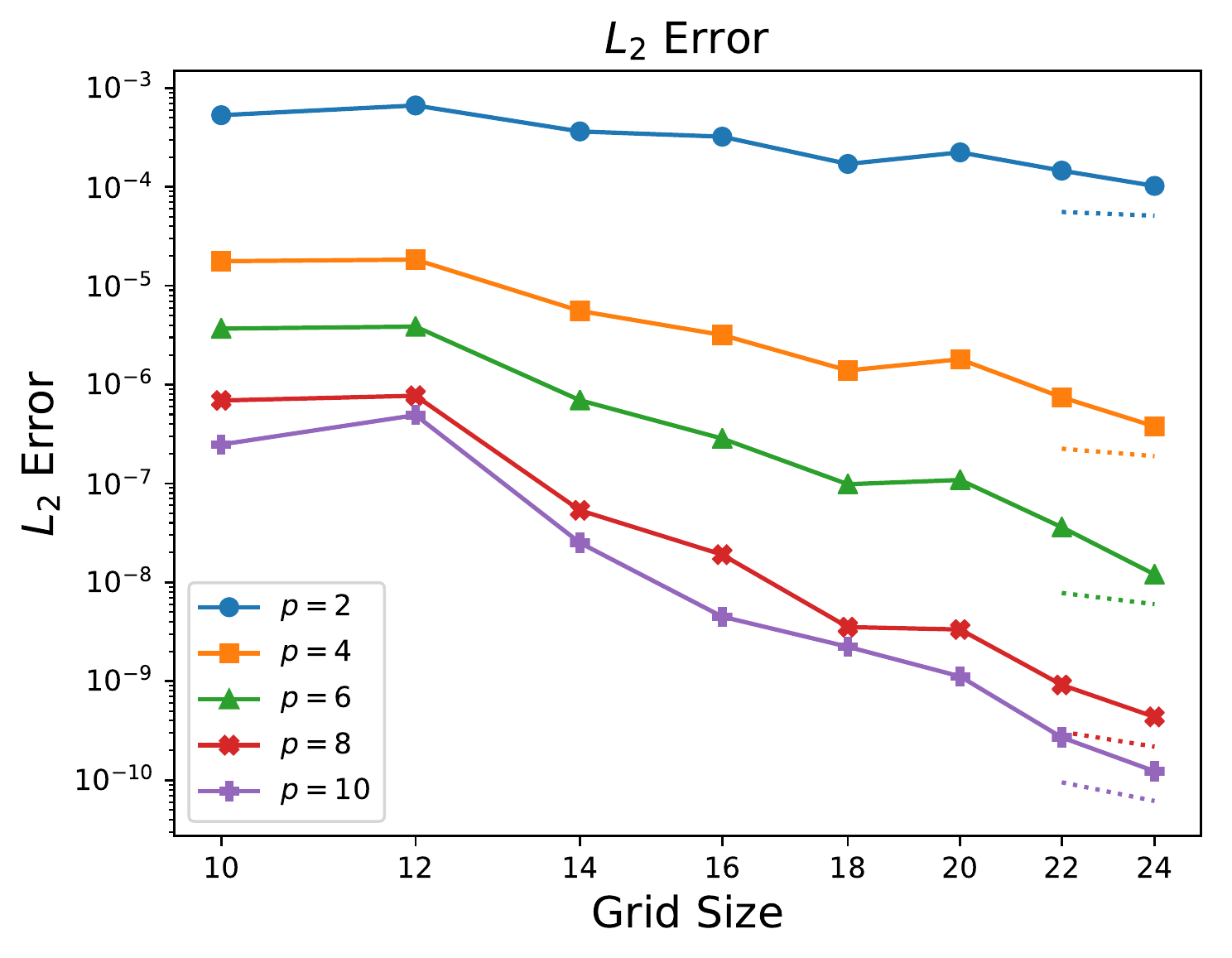}
\includegraphics[scale=.45]{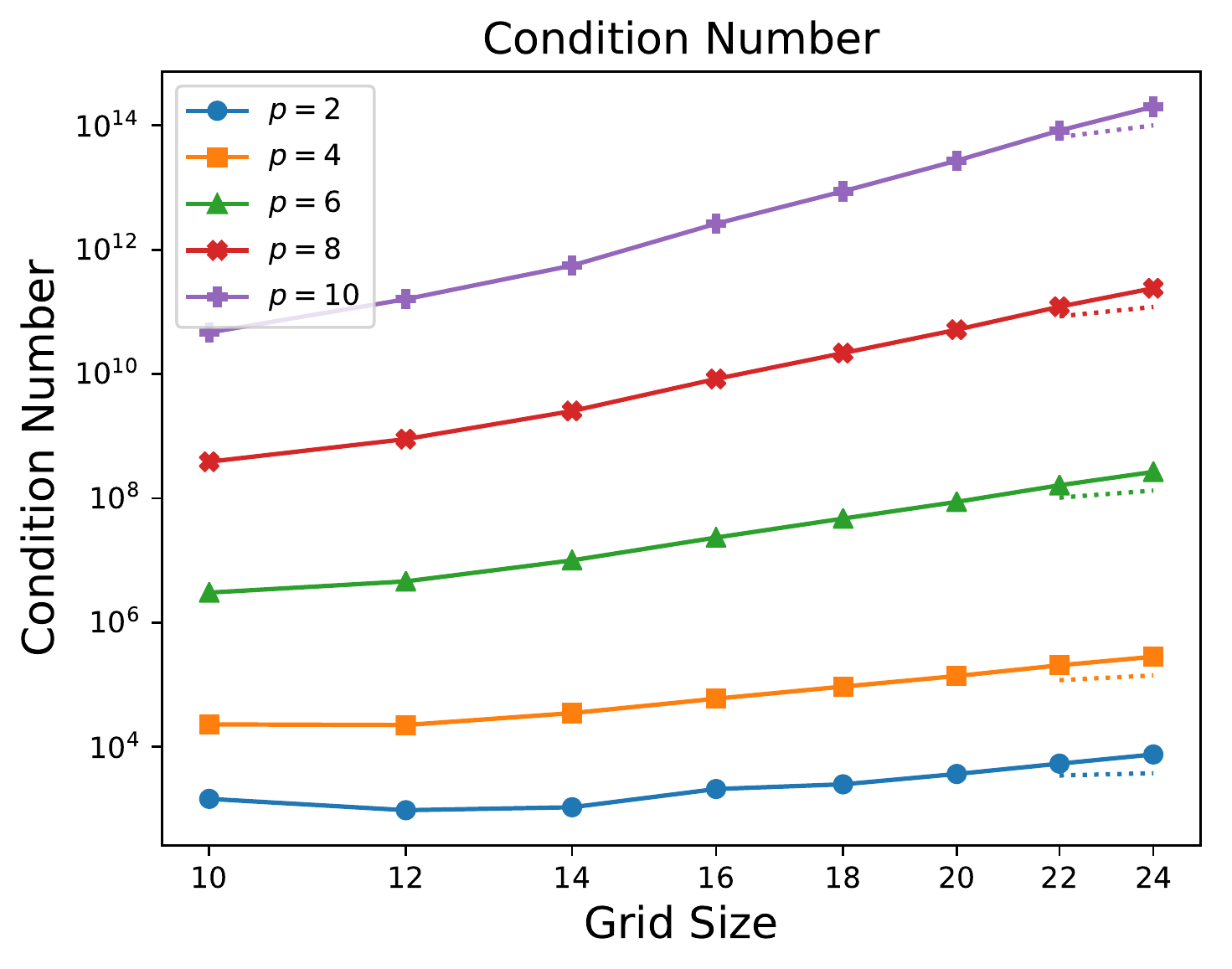}
\caption{Convergence of the $L_2$ error and growth of the condition number of the matrix $CS_p^{-1/2}$ for different values of $p$ when solving Equation \protect\eqref{eq:exppara}. The light dotted lines are reference lines of slope $\frac{1}{m^p}$ and $m^p$, where $m$ is the number of grid points along one dimension.}
\label{fig:ratespara}
\end{figure}

As in the elliptic case, the problem then reduces to finding 
$$
 u = S_p^{-1/2}(CS_p^{-1/2})^+ f.
$$
The solution is obtained using a QR decomposition of $S_p^{-1/2}C^\top$, as described in Section \ref{sec:soln}. The numerical results for the initial boundary value problem are summarized in Table \ref{tab:timespara} and Figure \ref{fig:ratespara}.
\begin{table}[H]
\centering
\begin{tabular}{|c|c|c|c|c|}
\hline
$|\mathbb{B}^m \times \mathbb{B}^n_E|$ & $|\Omega \times \widetilde{\mathbb{B}}^n_E|$ & $|\Omega^m \times \{0\}|$ & $|\Gamma^m \times \widetilde{\mathbb{B}}^n_E|$ &CPU Time\\
\hline 
$ 10 ^2 \times 11 = 1100$ & $ 280 $ & $ 28 $ & $ 130 $ & $ 1.2 $\\
\hline
$ 12 ^2 \times 11 = 1584$ & $ 400 $ & $ 40 $ & $ 150 $ & $ 1.9 $\\
\hline
$ 14 ^2 \times 11 = 2156$ & $ 500 $ & $ 50 $ & $ 180 $ & $ 3.6 $\\
\hline
$ 16 ^2 \times 11 = 2816$ & $ 720 $ & $ 72 $ & $ 200 $ & $ 4.7 $\\
\hline
$ 18 ^2 \times 11 = 3564$ & $ 860 $ & $ 86 $ & $ 230 $ & $ 7.2 $\\
\hline
$ 20 ^2 \times 11 = 4400$ & $ 1060 $ & $ 106 $ & $ 250 $ & $ 10.2 $\\
\hline
$ 22 ^2 \times 11 = 5324$ & $ 1340 $ & $ 134 $ & $ 280 $ & $ 18.7 $\\
\hline
$ 24 ^2 \times 11 = 6336$ & $ 1540 $ & $ 154 $ & $ 300 $ & $ 22.7 $\\
\hline
\end{tabular}
\caption{Grid sizes and computation times for solving Equation \protect\eqref{eq:exppara}. All computations were performed on an Intel I7-6660U.}
\label{tab:timespara}
\end{table}
\appendix
\numberwithin{equation}{subsection}
\section{The Chebyshev Polynomials}\label{sec:cheb}
Setting $\mathbb{B} = [-1,1]$, the Chebyshev polynomials of the first kind are given by 
$$
T_m(x) = \cos\bigl( m \arccos(x)\bigr),\: x \in \mathbb{B}
,\: m\in \mathbb{N}.
$$
We briefly describe some relevant properties of the Chebyshev polynomials used in the body of the paper. We refer to \cite{Tre00} for a more detailed discussion. 
\subsection{The Chebyshev roots}\label{sec:grid}
For fixed $m \in \mathbb{N}$, the $m$ roots of $T_m(x)$ are given by
$$
x_k = \cos \Bigl( \pi \frac{2k - 1}{2m} \Bigr),\: 0 \leq k \leq m-1.
$$
The Chebyshev grid comprising all roots of $T_m$, given by $\{ x_k\, |\, k=0,\dots,m-1\}$, is well adapted to the spectral calculation of derivatives.
In higher dimensions, a tensor product of one-dimensional Chebyshev grids
can be used. Throughout the body  of the paper, $\mathbb{B}^m$ has been used to denote the Chebyshev grid of the appropriate dimension.
\subsection{Orthogonality relations}
The sequence $(T_m)_{m \in \mathbb{N}}$ forms an orthogonal basis for $L_2(\mathbb{B})$ with respect to the measure $\frac{dx}{\sqrt{1-x^2}}$. More specifically, for $i,j \in \mathbb{N}$, 
$$ \int_{-1}^1 T_i(x)T_j(x) \frac{dx}{\sqrt{1-x^2}} = 
\begin{cases}
0,&\text{if } i \neq j,\\
\pi,&\text{if } i = j = 0,\\
\pi/2,&\text{if } i = j \neq 0.
\end{cases}$$
The Chebyshev functions restricted to $\mathbb{B}^m$ also satisfy a discrete orthogonality relation.
Indeed, for $0 \leq i,j \leq m-1$, one has that
$$\sum_{k=0}^{m-1} T_i(x_k) T_j(x_k) = 
\begin{cases}
0,&\text{if } i \neq j, \\
m,&\text{if } i = j = 0, \\
\frac{m}{2},&\text{if } i = j \neq 0.
\end{cases}
$$
\subsection{The Chebyshev transform}\label{sec:transform}
Because $(T_m)_{m\in \mathbb{N}}$ forms an orthogonal basis of
$\operatorname{L} _2(\mathbb{B})$,
any function $u\in \operatorname{L} _2(\mathbb{B})$ can be developed in a ``Chebyshev series''.
We set
$$c_k = \frac{p_k}{\pi} \int_{-1}^1 u(x) T_k(x) \frac{dx}{\sqrt{1-x^2}}, 
\text{ where }
p_k = \begin{cases}
1,&\text{if } k = 0, \\
2,&\text{if } k \neq 0.
\end{cases}
$$
so that
$$ u(x) = \sum_{m=0}^{\infty} c_k T_k(x).$$
It is also possible to define the Chebyshev transform, denoted by $\mathfrak{C}$, which maps a function to the sequence of its Chebyshev coefficients.
$$
 \mathfrak{C}(u) = (c_k)_{k\in \mathbb{N}}.
$$
The discrete orthogonality relation also yields  a discrete version of a Chebyshev expansion. Given $u:\mathbb{B}^m \rightarrow \mathbb{R}$, let
$$ c_k = \frac{p_k}{m}\sum_{i=0}^{m-1} u_i T_k(x_i)
\text{ for }
p_k = \begin{cases}
1,&\text{if } k = 0,\\
2,&\text{if } k \neq 0.
\end{cases}
$$
Then $u$ can be written as a discrete Chebyshev series
$$
u_{\bullet} = \sum_{k=0}^{m-1} c_k T_k(x_{\bullet}).
$$
As in the continuous case, we can define the discrete Chebyshev transform $\mathfrak{C}_m$,
which maps $u$ to its discrete Chebyshev series, i.e., we set
$$
\mathfrak{C}_m(u) = (c_k)_{k=0,\dots,m-1}=:c.
$$
The discrete Chebyshev transform, $\mathfrak{C}_m$, as well as its inverse $\mathfrak{C}_m^{-1}$, can be implemented efficiently using an FFT algorithm, in the form of the discrete cosine transform. More specifically,
$$\mathfrak{C_m}(u)_k = a_k\textsf{DCT}(u)_k \text{ for } a_k = 
\begin{cases} 
\frac{1}{2m},&\text{if } k = 0, \\
\frac{1}{m},&\text{if } k > 0,
\end{cases}$$
and
$$
 \mathfrak{C_m}^{-1}(c) =
 \textsf{IDCT}(\,\widetilde{c}\,)\text{ for } \widetilde{c}_k= 
\begin{cases} 
c_k,&\text{ if } k = 0, \\
c_k/2,&\text{ if } k > 0.
\end{cases}
$$
In dimension larger than one, $\mathfrak{C}$ and $\mathfrak{C}_m$
will denote the continuous and discrete one-variable Chebyshev transforms applied
successively in each direction. Numerically, this can be accomplished with the use of
$\textsf{DCTN}$, where the factors $a_k$ and $b_k$ are raised to
the power of the dimension.

\subsection{Derivative formulas}\label{sec:derivatives}
Discrete derivatives can be efficiently evaluated on the Chebyshev grid using the $\textsf{DCT}$ and $\textsf{DST}$. We denote by $(k_{\bullet})_{k_{\bullet} \in \{0, \ldots, m-1\}}$ the discrete frequency
vector, and let $M[\bullet]$ represent multiplication by the discrete function $\bullet$. We also define a shifting operator $\textsf{R}$ with 
$$
\textsf{R}_{ij} = \delta_{i+1,j},
$$
so that $\textsf{R}$ is the matrix with ones on the superdiagonal.
Then, given a function $u=(u_i)_{i \in \{0,\dots,m-1 \}}$ defined on the Chebyshev grid $\mathbb{B}^m$,
a spectrally accurate discrete derivative $Du$ can be calculated using the matrix given by
$$
  D =M \Big[ \frac{1}{\sqrt{1-x_{\bullet}^2}} \Big] \circ \textsf{IDST}\circ\textsf{R} \circ M \Big[ \frac{k_{\bullet}}{2m} \Big] \circ\textsf{DCT}.
$$
Similarly, we can compute
$$ D^2 = M \Big[ \frac{-1}{1-{x_{\bullet}}^2} \Big] \circ \textsf{IDCT} \circ M \Big[ -\frac{k_{\bullet}^2}{2m} \Big] \circ \textsf{DCT} + M \Big[ \frac{x_{\bullet}}{(1-x_{\bullet}^2)^{3/2}} \Big] \circ \textsf{IDST} \circ \textsf{R} \circ M \Big[ -\frac{k_{\bullet}}{2m}\Big] \circ \textsf{DCT}. $$
Of course, a corresponding operator can be formed in higher dimensions, where the $\textsf{DCT}$, $\textsf{DST}$ as well as the frequency vector $k_{\bullet}$ are taken along the desired direction of differentiation. In the body of the paper, the derivative operator in the $x_i$ direction was denoted as $D_i$.

\subsection{Eigenvalue equation}\label{sec:eig}

Setting $\mathcal{D} = M \bigl[ \sqrt{1-x^2} \bigr] \circ \frac{\partial}{\partial x} $, the Chebyshev polynomials satisfy the eigenvalue equation
\begin{equation}\label{eq:eig}
-\mathcal{D}^2 T_m(x) = m^2 T_m(x).
\end{equation}
Similarly, given the Chebyshev grid $\mathbb{B}^m$, the discrete Chebyshev functions $T_j(x_\bullet)$ satisfy a discrete eigenvalue equation. Defining
$\mathcal{D}_m = M \bigl[ \sqrt{1-x_{\bullet}^2}\bigr] \circ D$, $T_j(x_\bullet)$ satisfies
$$
-\mathcal{D}_m^2 T_j(x_\bullet) = j^2 T_j(x_\bullet),\: j \in \{ 0, \ldots, m-1 \}.
$$
This implies that
$$
(1-\mathcal{D}_m^2)^{-p/2}= \mathfrak{C}_m^{-1} \circ M \Big[(1+|k_{\bullet}|^2)^{-p/2} \Big] \circ\mathfrak{C}_m.
$$

\subsection{Interpolation operators}\label{sec:interpolation}

Functions defined on the Chebyshev grid can be interpolated at arbitrary points in $\mathbb{B}$. Such interpolation can be stably computed by means of the barycentric interpolation
formulas described in \cite{Tre04}. Define first the vector $w$ by
$$
w_k = (-1)^k \sin \Bigl( \frac{2k-1}{2m} \Bigr),\:0 \leq k \leq m-1.
$$
If $y \in \mathbb{B}$ and $(x_i)_{i \in \{0,\dots,m-1\}}$ is the vector of points in $\mathbb{B}^m$, a spectrally accurate interpolation of a discrete function $u$ defined on the Chebyshev grid $\mathbb{B}^m$ can be obtained by
$$
u(y) = \delta_y \cdot u \text{ where } (\delta_y)_i = \frac{1}{\sum_{k=0}^{m-1} \frac{w_k}{y-x_k}} \frac{w_i}{y-x_i}.
$$
To calculate the interpolation of the first derivative, which will be used in the Neumann problem, we use the derivative of the above formula,
$$
Du(y) = \delta_{y} \circ D \text{ where } (\delta_{y} \circ D)_i  = - \frac{1}{\sum_{k=0}^{m-1} \frac{w_k}{y-x_k}} \frac{w_i}{(y-x_i)^2} + 
\frac{\left( \sum_{k=0}^{m-1} \frac{w_k}{(y-x_k)^2} \right)}
{\left( \sum_{k=0}^{m-1} \frac{w_k}{y-x_k} \right)^2}
\frac{w_i}{y-x_i}.$$
To interpolate in several dimensions, we use a tensor product of the given interpolations, which were denoted $\delta_y$ and $(\delta_y \cdot \nabla)$. To calculate a directional derivative of the grid function
$u$ in the direction $\nu$ at the point $y$, we use $(\delta_y \circ \nabla u) \cdot \nu_y$.

\bibliography{cheb}

\begin{thebibliography}{1}

\bibitem{AG18}
D~Agress and P~Guidotti.
\newblock A novel optimization approach to fictitious domain methods.
\newblock {\em arXiv preprint arXiv:1808.02158}, 2018.

\bibitem{Tre04}
JP~Berrut and LN~Trefethen.
\newblock Barycentric lagrange interpolation.
\newblock {\em SIAM review}, 46(3):501--517, 2004.

\bibitem{Fas05}
GE~Fasshauer.
\newblock Meshfree methods.
\newblock {\em Handbook of theoretical and computational nanotechnology},
  27:33--97, 2005.

\bibitem{Fas07}
GE~Fasshauer.
\newblock {\em Meshfree approximation methods with MATLAB}, volume~6.
\newblock World Scientific, 2007.

\bibitem{Glo94}
R~Glowinski, TW~Pan, and J~Periaux.
\newblock A fictitious domain method for dirichlet problem and applications.
\newblock {\em Computer Methods in Applied Mechanics and Engineering},
  111(3-4):283--303, 1994.

\bibitem{Gon10}
{\'A}lvaro Gonz{\'a}lez.
\newblock Measurement of areas on a sphere using fibonacci and
  latitude--longitude lattices.
\newblock {\em Mathematical Geosciences}, 42(1):49, 2010.

\bibitem{Pal16}
Richard Palais, Bob Palais, and Hermann Karcher.
\newblock Pointclouds: Distributing points uniformly on a surface.
\newblock {\em arXiv preprint arXiv:1611.04690}, 2016.

\bibitem{Sch07}
R~Schaback.
\newblock Kernel-based meshless methods.
\newblock {\em Lecture Notes for Taught Course in Approximation Theory.
  Georg-August-Universit{\"a}t G{\"o}ttingen}, 2007.

\bibitem{Tre00}
LN~Trefethen.
\newblock {\em Spectral methods in MATLAB}, volume~10.
\newblock Siam, 2000.

\end{thebibliography}
\end{document}